%% file: HIEvsPML.tex
\documentclass{elsarticle}

\usepackage{amsfonts}
\usepackage{amssymb,amsmath} 
\usepackage{graphicx} 
\usepackage{subfigure}
\usepackage{color}
\usepackage{cellspace} 
\usepackage{longtable}

\newtheorem{Theo}{Theorem}
\newdefinition{Def}[Theo]{Definition}
\newtheorem{Lem}[Theo]{Lemma}
\newdefinition{Rem}[Theo]{Remark}

\newcommand{\setR}{\mathbb{R}} 
\newcommand{\setC}{\mathbb{C}}

\newcommand{\calL}{{\cal L}} 
\newcommand{\calM}{{\cal M}}
\newcommand{\calT}{{\cal T}} 
\newcommand{\calD}{{\cal D}}

\newcommand{\bfg}{{\bf g}}

\newcommand{\Oi}{\Omega_{\rm int}} 
\newcommand{\Oe}{\Omega_{\rm ext}}

\newcommand{\MT}{\calM_{\kappa_0}} 
\newcommand{\OpT}{\calT}
\newcommand{\BilF}[2]{A\!\left(#1,#2\right)}

\newcommand{\Dop}{\calD} \def\id{{\rm id}}
\def\vector#1#2{\left(\!\!\begin{array}{c} #1 \\ #2 \end{array}\!\!\right)}

\begin{document}

\begin{frontmatter}
\title{Hardy space infinite elements for Helmholtz-type problems with unbounded inhomogeneities \tnoteref{label1} }
\tnotetext[label1]{This work
  was supported by the Deutsche Forschungsgemeinschaft (DFG).}
\author[ln]{Lothar Nannen}
\author[as]{Achim Sch\"adle}
\address[ln]{Institut f\"ur Numerische und Angewandte Mathematik,
Georg-August Universit\"at G\"ottingen,
Lotzestra\ss e 16-18,
37083 G\"ottingen, Germany}
\address[as]{Mathematisches Institut, Heinrich-Heine Universit\"at D\"usseldorf,Universit\"atsstra\ss e 1, 
40225 D\"usseldorf, Germany}

\begin{abstract} 
This paper introduces a class of approximate transparent boundary
conditions for the solution of Helmholtz-type resonance and scattering
problems on unbounded domains. The computational domain is assumed to
be a polygon. A detailed description of two variants of the Hardy space infinite
element method which relays on the pole condition is given. The method
can treat waveguide-type inhomogeneities in the domain with
non-compact support. The results of the Hardy space infinite element
method are compared to a perfectly matched layer method. Numerical
experiments indicate that the approximation error of the Hardy space
decays exponentially in the number of Hardy space modes.
\end{abstract}

\end{frontmatter}

\section{Introduction} 
\label{sec:intro}
To solve numerically the Helmholtz equation 
\begin{equation}
  \label{eq:inTHelmholtz}
  - \Delta u(x) - \kappa^2~n(x) u(x) = 0,\qquad x\in \Omega
\end{equation}
on an unbounded domain $\Omega$ with some boundary conditions on
$\partial \Omega$ and a radiation condition at infinity, the
computational domain is typically restricted to a bounded interior
$\Omega_{\rm int}$. Here $n(x)$ is the refraction index and $\kappa$
is the wavenumber. In case of a scattering problem $\kappa>0$ is
given, whereas in case of a resonance problem $\kappa$ with positive
real part is the sought resonance. Applying transparent boundary
conditions at the artificial interface $\Gamma$ of the
exterior/interior domain, the problem can be restricted to solving the
Helmholtz equation on the interior domain only. Transparent boundary
conditions on $\Gamma$ have to model the correct radiation condition
at infinity.

For homogeneous exterior domains and $\kappa>0$ the correct radiation
condition is the Sommerfeld radiation condition. It is known that the
exact transparent boundary condition, i.e. the Calderon or
Dirichlet-to-Neumann map, is non-local. Currently used methods to
approximate or realize transparent boundary conditions are based on
separable coordinates and special functions~\cite{GroteKeller:98},
infinite elements~\cite{Astley:00,DemkowiczGerdes:98}, perfectly
matched layer (PML) constructions~\cite{Simon:79,Berenger:94,ChewWeedon:94},
boundary integral approaches~\cite{HsiaoWendland:08} and local high
order approximations~\cite{Givoli:04}.

Except for some PML methods they depend non-linearly on $\kappa^2$,
which is a severe drawback when solving resonance problems, where one
seeks non-trivial eigenpairs $(u,\kappa^2)$ for
\eqref{eq:inTHelmholtz} with vanishing boundary conditions on
$\partial \Omega$ and a radiation condition at
infinity. Using finite elements to discretize~\eqref{eq:inTHelmholtz}
in $H^1(\Omega_{\rm int})$ leads to a generalized eigenvalue problem.
Employing transparent boundary conditions, that are non-linear in the
eigenvalue, would make the eigenvalue problem non-linear. Although it
is possible to solve the resulting non-linear eigenvalue problems, see e.g.
\cite{LenoirVulliermeHazard:92}, it is reasonable to avoid them.
Therefore PML methods are currently the standard
method for solving resonance problems, see
e.g.~\cite{HeinHohageKochSchoeberl:07,KimPasciak:08}. Under the name
{\it complex scaling} they have been used since the 1970s for the
theoretical study and the numerical computation of resonances in
molecular physics \cite{Simon:79}. Unfortunately these methods give
rise to spurious resonance modes and several parameters have to be
optimized for each problem. The PML method we use
to check the results of our Hardy space infinite element method choses
the thickness and the discretization of the layer
adaptively~\cite{ZschiedrichBurgeretal:05,SchaedleZschiedrichetal:06}.

The theoretical framework of the Hardy space infinite element (HSIE) method is the
pole condition by F. Schmidt~\cite{Schmidt:98,Schmidt:02}. The pole condition considers 
the Laplacetransform of the exterior solution with respect to some 
generalized distance variable. A solution is then called (purely) outgoing if this Laplacetransform 
has no singularity in the lower complex half plane, vice versa a solution is (purely) incoming if 
its Laplacetransform has no singularity in the upper complex half plane. 
In~\cite{HohageSchmidtZschiedrich:03a} it is shown that for
homogeneous exterior domains this condition on the singularities of the
Laplacetransform, i.e. all its singularities are located in the upper complex half, 
is equivalent to the Sommerfeld radiation condition.

Compared to former numerical realizations of the real axis approach of the pole
condition~\cite{Hohageetal:02}, which were based on BDF and Runge-Kutta
methods, the HSIE method which is based on a Galerkin method in the Hardy space 
$H^+(D)$ of the complex unit disk, shows
exponential convergence and is simple to implement in finite element codes. It was
first presented in \cite{Nannen:08,HohageNannen:08} for homogeneous exterior domains
with spherical interface $\Gamma$.  The cut function approach of the pole
condition~\cite{Hohageetal:02} was discretized in \cite{Schmidtetal:07} using a
collocation method in some sort of Hardy space and shows exponential convergence in
experiments, too.  Moreover it allows for the evaluation of the exterior field. 
However this approach is not linear in $\kappa^2$ and exhibits problems concerning the
stability under perturbations of the boundary data.

Extending~\cite{HohageNannen:08} we assume here that the interface $\Gamma$ is the 
boundary of a convex polygon $P$, i.e. $\Omega_{\mathrm{int}} = P \cap \Omega$ 
and $\Omega_{\mathrm{ext}} = \setR^2 \setminus P$, and that the exterior domain 
$\Omega_{\mathrm{ext}}$ is discretized by
infinite trapezoids, such that the refractive index $n(x)$ 
of~\eqref{eq:inTHelmholtz} is constant on each trapezoid. These trapezoids are
the images of bilinear mappings from a reference strip. The infinite
direction of this reference strip is mapped into the Hardy space $H^+(D)$, where a
$L^2$-orthogonal basis is given by trigonometric monomials. The basis functions of
our Galerkin method are therefore tensor products of standard finite element
functions with the trigonometric monomials in the Hardy space.

The method presented here treats scattering problems as well as the corresponding
resonance problems. The pole condition as a mean to realize transparent boundary 
conditions for time dependent problems is considered in~\cite{Ruprechtetal:08}. The
discretization used there is almost equivalent to one employed here. 

In section 2 we present the HSIE method from a practical
point of view first in one dimension and then in two dimensions. In section 3 we
shortly review the PML method used here. Numerical results are
presented in section 4 with a comparison of both methods.

\section{Hardy space infinite element method}
\label{sec:HSIEM}
To explain the basic ideas of HSIEs we shortly discuss
one-dimensional problems, even though for such problems there exist simpler and more efficient
methods to treat the unboundedness of the domain. Nevertheless, as the
multi-dimensional elements are tensor products of standard finite elements and
one-dimensional infinite elements, it is useful to start with the
simple case.

\subsection{One-dimensional elements}
\label{sec:HSIEMone}
We consider the Helmholtz equation
\begin{subequations}
  \label{pr:1dim}
  \begin{eqnarray}
    -u''(r) - \kappa^2 n(r) u(r) &=& 0,\quad r \geq 0,
    \label{eq:Helmholtz1d} 
    \\ 
    u'(0) &=&g,
    \\ 
    u &&\text{is outgoing}
  \end{eqnarray}
\end{subequations} 
with complex wave number $\kappa \in \setC$ with positive real part ($\Re(\kappa)>0)$, 
boundary value $g\in\setC$, and positive potential $n \in L^\infty((0,\infty))$ satisfying
$n(r)=1$ for $r\geq a$ for some $a$. The Sommerfeld radiation condition
\begin{equation*}
\lim_{r \to \infty} \left(\partial_{r} u(r) -i\kappa u(r)\right) =0\label{eq:Sommerfeld1d}
\end{equation*}
guarantees for $\kappa$ with nonnegative imaginary part ($\Im(\kappa) \ge 0$ 
that~\eqref{pr:1dim} is well-posed and that the solution $u$ is outward
radiating. Solutions to \eqref{eq:Helmholtz1d} may be decomposed into an interior part 
$u_{\rm int}:=u|_{[0,a]}$ and an exterior part $u_{\rm ext}(r):= u(r+a)$, $r\geq0$, with
\begin{equation*} 
  u_{\rm ext}(r)=C_1 e^{i \kappa r} + C_2 e^{-i \kappa r} \quad
  \text{and}~C_1+C_2=u_{\rm int}(a).
\end{equation*}
The term $C_1 e^{i \kappa r}$ corresponds to an outgoing radiating wave, that satisfies
the Sommerfeld radiation condition, and $C_2 e^{-i \kappa r}$ corresponds to an incoming wave. 
Therefore $C_{2}=0$ and hence $u_{\rm ext}(r)=u_{\rm int}(a) e^{i \kappa r}$.
The solution of~\eqref{pr:1dim} restricted to $[0,a]$ is thus given by the 
simple boundary value problem
\begin{equation}
  \label{pr:1dimint} 
  -u_{\rm int}''(r) - n(r) \kappa^2 u_{\rm int}(r) = 0,\qquad u_{\rm
    int}'(0)=g, \qquad u_{\rm int}'(a) =i\kappa u_{\rm int}(a).
\end{equation} 
Thus in the one-dimensional case the Sommerfeld radiation condition 
is satisfied at the boundary of the interior domain. Therefore it can be used for $\kappa$ with negative imaginary part, too.

Note that the resonance problem corresponding
to \eqref{pr:1dimint} leads to a quadratic eigenvalue problem in $\kappa$ 
whereas the discretization using Hardy space infinite elements will lead to 
an eigenvalue problem that is linear in $\kappa^2$.

Hardy space infinite elements rest on the fact, 
that the Laplace transform of the exterior solution
\begin{equation*} 
  \left(\calL u_{\rm ext}\right) (s) :=\int_0^\infty e^{-sr}
  u_{\rm ext}(r) dr,\qquad \Re(s) \geq \left|\Im(\kappa)\right|
\end{equation*} 
has a holomorphic extension except for two poles at $\pm i \kappa$
\begin{equation*} 
  \calL \left\{C_1 e^{i \kappa \bullet} + C_2 e^{-i \kappa \bullet}\right\}(s) =
  \frac{C_1}{s - i \kappa} + \frac{C_2}{s + i \kappa}.
\end{equation*} 
Hence, $u_{\rm ext}$ is outgoing if and only if $\calL u_{\rm ext}$ has no poles 
with negative imaginary part. An equivalent formulation in terms of an appropriate
Hardy space is given in Definition~\ref{Def:Pole}.

\begin{figure}
  \label{fig:Hardy} \centering
  \resizebox{0.8\textwidth}{!}{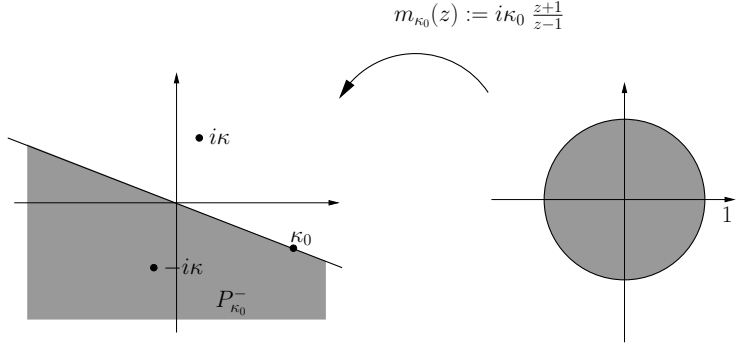}
  \caption{sketch to definition \ref{Def:Hardy}}
\end{figure}

\begin{Def}[Hardy space]
  \label{Def:Hardy} 
  Let $P_{\kappa_0}^-=\{ s \in \setC : \Im(s/\kappa_0)<0\}$ be the half plane below the
  line $\kappa_0 \setR$ through the origin and $\kappa_{0}$, see Fig. \ref{fig:Hardy}. The \emph{Hardy
    space} $H^-(P_{\kappa_{0}}^-)$ is the space of all functions $f$, that are holomorphic
  in $P_{\kappa_0}$, such that 
  $$\int_{\setR} |f(\kappa_0x- \kappa_0i \epsilon)| ^2 dx$$ is
  uniformly bounded for $\epsilon>0$.
  \\
  Let $D =\{ s \in \setC : |s|<1\} $ be the open unit disk. The \emph{Hardy space}
  $H^{+}(D)$ is the space of all functions $f$, that are holomorphic in $D$, such that
  $$\int_{0}^{2\pi} |f(r \mathrm{e}^{it})|^{2} dt$$
  is bounded uniformly for $r \in [0,1)$.
\end{Def}

Due to the uniform boundedness of $f$ there exist in both cases a
$L^2$ function on the boundary, which is uniquely determined by $f$ and which
determines vice versa uniquely the function in the domain. Hence, we
identify a Hardy space function $f \in H^-(P_{\kappa_{0}}^-)$ or $f
\in H^+(D)$ with its boundary function $f \in L^2(\kappa_0\setR)$ and
$f \in L^2(S^1)$ respectively.

\begin{Def}[Pole condition]
  \label{Def:Pole} Let $\kappa_0$ be a complex constant with positive real part and
  $\Re(\kappa/\kappa_0)>0$. Then a solution $u$ to \eqref{eq:Helmholtz1d} is said to
  obey the pole condition and is called outgoing, if the holomorphic extension of the
  Laplace transform of the exterior part lies in the Hardy space $H^-(P_{\kappa_0}^-)$.
\end{Def}

The constant $\kappa_0$ will act as a tuning parameter in the
method presented below. In future we omit the formulation ''holomorphic extension 
of the Laplace transform'' and shortly write Laplace
transform $\calL$. The Hardy space $H^-(P_{\kappa_0}^-)$ is a Hilbert space
\cite{Duren:70,Hoffman:62}, with the standard $L^2$-norm, and the following Lemma connects the two 
Hardy spaces $H^-(P_{\kappa_{0}}^-) $ and $ H^+(D)$.

\begin{Lem}[M\"obius transform]
  \label{Lem:Moebius} The mapping
  \begin{equation} 
   \MT\ : \ H^-(P_{\kappa_{0}}^-) \to H^+(D) \ : \
    f \mapsto (\MT f)(z):= f\left( i \kappa_0
      \frac{z+1}{z-1}\right) \frac{1}{z-1}
  \end{equation} 
  is  up to a factor $\sqrt{2 | \kappa_0|}$ unitary. 
\end{Lem} 

Due to the explicit knowledge of $u_{\rm ext}$ the transformed
function $\hat{U}:= \MT \calL u_{\rm ext}$ is given by
\begin{equation}
  \label{eq:reprUhat} 
  \hat{U}(z)=\frac{u_{\rm int}(a)}{i \kappa_0 (z+1) -
    i \kappa (z-1)}=\frac{u_{\rm int}(a)}{i(\kappa+\kappa_0)} \sum_{j=0}^\infty
  \left(\frac{\kappa-\kappa_0}{\kappa+\kappa_0} \right)^{j} z^j.
\end{equation}
Since $\left|\frac{\kappa-\kappa_0}{\kappa+\kappa_0} \right|<1$, we
could expect exponential convergence for the exterior solution, if we
use the first $N+1$ trigonometric monomials $\{z^0,z^1,...,z^N\}$ as a
Galerkin basis of the space $H^+(D)$. For the interior part $u_{\rm
int} \in H^1([0,a])$ we use a standard finite element method.

Both methods have the term $u_{\rm int}(a)$ in common and we call the
associated boundary degree of freedom $u_0$. As usual it couples the
finite elements for the interior part with the infinite elements for
the exterior part. Note, that in the formulation \eqref{eq:reprUhat}
all degrees of freedom for $\hat{U}$ would couple with $u_{\rm
int}(a)$, since $u_{\rm int}(a)=2 i \kappa_0 \hat{U}(1)$. In order to
get a local coupling of the boundary degree of freedom and the inner
degrees of freedom, we decompose $\hat{U}=\frac{1}{ i
\kappa_0}\OpT_-(u_0 ,U)^\top$ with
\begin{equation}
  \label{eq:defOpT} \OpT_-\vector{u_0}{U}(z):=\frac{1}{2}\left(u_0 + (z-1)
    U(z) \right)
\end{equation} 
for $(u_0,U)^\top \in \setC \times H^+(D)$ and use the trigonometric
monomials as a basis for $U$.

It remains to derive a variational
formulation for $(u_{\rm int},U)^\top$ in $H^1([0,a]) \times H^+(D)$,
which was done in~\cite{HohageNannen:08}. Here $H^1([0,a])$ is the
Sobolev space of weakly differentiable functions on the interval
$[0,a]$.  The basic idea is to use the properties of the Fourier
transform to get the identity
\begin{equation}
  \label{eq:basisid} 
  \int_0^\infty f(r) g(r) dr = -2 i \kappa_0 \BilF{\MT \calL
    f}{\MT \calL g}
\end{equation} 
with
\begin{equation} 
  \label{eq:BiLF}
  \BilF{F}{G}:=\frac{1}{2 \pi}\int_{S^1} F(\overline{z}) G(z)
  |dz|,\qquad F,G \in H^+(D).
\end{equation} 
This holds for $u_{\rm ext}$ and suitable test functions 
$v_{\rm ext}$, as well as for the derivatives $u_{\rm ext}'$ and $v_{\rm ext}'$. 
Using the decomposition $\MT \calL u_{\rm ext}=\frac{1}{i \kappa_0} \OpT_-(u_0,U)^\top$
and $\MT \calL v_{\rm ext}=\frac{1}{i \kappa_0} \OpT_-(v_0,V)^\top$ we obtain simple
formulas for the derivatives of $u_{\rm ext}$ and $v_{\rm ext}$
\begin{equation}
  \label{eq:OpTp} 
  \MT \calL f'=\OpT_+\vector{f_0}{F} \quad \text{with}~
  \OpT_+\vector{f_0}{F}(z):=\frac{1}{2}\left(f_0 + (z+1) F(z) \right).
\end{equation} 
Now we are able to deduce from the formal variational formulation
of \eqref{eq:Helmholtz1d}
\begin{equation*} \int_0^a (u_{\rm int}' v_{\rm int}' - \kappa^2 n\, u_{\rm int}
  v_{\rm int}) dr + \int_0^\infty (u_{\rm ext}'v_{\rm ext}'-\kappa^2 n u_{\rm ext}
  v_{\rm ext}) dr = - g v_{\rm int}(0)
\end{equation*} 
and \eqref{eq:basisid} the variational equation in $H^1([0,a])
\times H^+(D)$:
\begin{equation}
  \label{1dimVariation} 
  B\left(\vector{u_{\rm int}}{U},\vector{v_{\rm int}}{V}\right) = -gv_{\rm int}(0)
\end{equation} with
\begin{align*} 
  & B\left(\vector{u_{\rm int}}{U},\vector{v_{\rm int}}{V}\right)
  := \int_0^a (u_{\rm int}' v_{\rm int}'-\kappa^2n\, u_{\rm int} v_{\rm int})dr \\
  &\quad - 2 i\kappa_0 \BilF {\OpT_+ \vector{u_0}{U}}{\OpT_+ \vector{v_0}{V}}
  -\frac{2 i\kappa^2}{\kappa_0}\BilF {\OpT_- \vector{u_0}{U}}{\OpT_-
    \vector{v_0}{V}}.
\end{align*} 
Since for the trigonometric monomials
$\BilF{z^j}{z^k}=\delta_{j,k}$, the implementation of the exterior part of the
bilinear form $B$ reduces to the implementation of the operators
$\OpT_\pm:\setC \times H^+(D) \to H^+(D)$, if the finite dimensional ansatz
space $\Pi_N:={\rm span}\{z^0,z^1,...,z^N\}$ is used for $H^+(D)$:
\begin{equation}
  \label{eq:OpTN} 
  \OpT_{N,\pm}:=\frac{1}{2}{\scriptstyle
    \left(\begin{array}{cccc} 1 & \pm1 & & \\ & \ddots & \ddots &\\ & & 1 & \pm1\\ &
        & & 1
      \end{array}\right)}.
\end{equation} 
The first row in these matrices correspond to the boundary
degree of freedom $u_0$. The local element matrix for the infinite element is
then given by
\begin{equation*} 
  -2i\kappa_0\OpT_{N,+}^{\top}\OpT_{N,+} -\kappa^2 \frac{2
    i}{\kappa_0}\OpT_{N,-}^{\top}\OpT_{N,-}.
\end{equation*} 
Note that it is linear in $\kappa^2$ and the matrices are
explicitly known. In \cite{HohageNannen:08} the equivalence of the variational equation in
$H^1([0,a])\times H^+(D)$ and the classical problem \eqref{pr:1dim} is
shown. Moreover, the stability of the Hardy space infinite element method
follows with a G\aa rding inequality and exponential convergence in the number
of degrees of freedom $N$ for the Hardy space is proven.
\begin{Rem} 
  In the space domain the monomial basis functions $z^j$ correspond to
  the functions
  \begin{equation*} 
    u_j(r)=e^{i \kappa_0 r} \left\{ u_0 + \sum_{n=0}^j
      \left(\begin{array}{c} j\\ n\end{array}\right) \frac{ (2 i \kappa_0
        r)^{n+1}}{(n+1)!} \right\}.
  \end{equation*} 
  Therefore for the one dimensional problem an optimal choice for scattering problems is
  $\kappa_0=\kappa$ since in this case the exact transparent boundary condition is
  obtained even with no degrees of freedom in $H^+(D)$. For resonance problems,
  $\kappa_0$ should be chosen in the region of the complex plane where resonances
  are of interest.
\end{Rem}

\subsection{Tensor product elements for a semi-infinite strip}
\label{sec:TensorStrip}

\begin{figure}
\begin{center}
  \subfigure[\label{fig:Strip} ]{\includegraphics[width=0.67\textwidth]{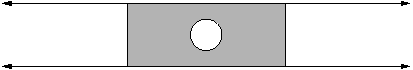}} 
  \subfigure[ \label{fig:WGpoles}]{\resizebox{.3\textwidth}{!}{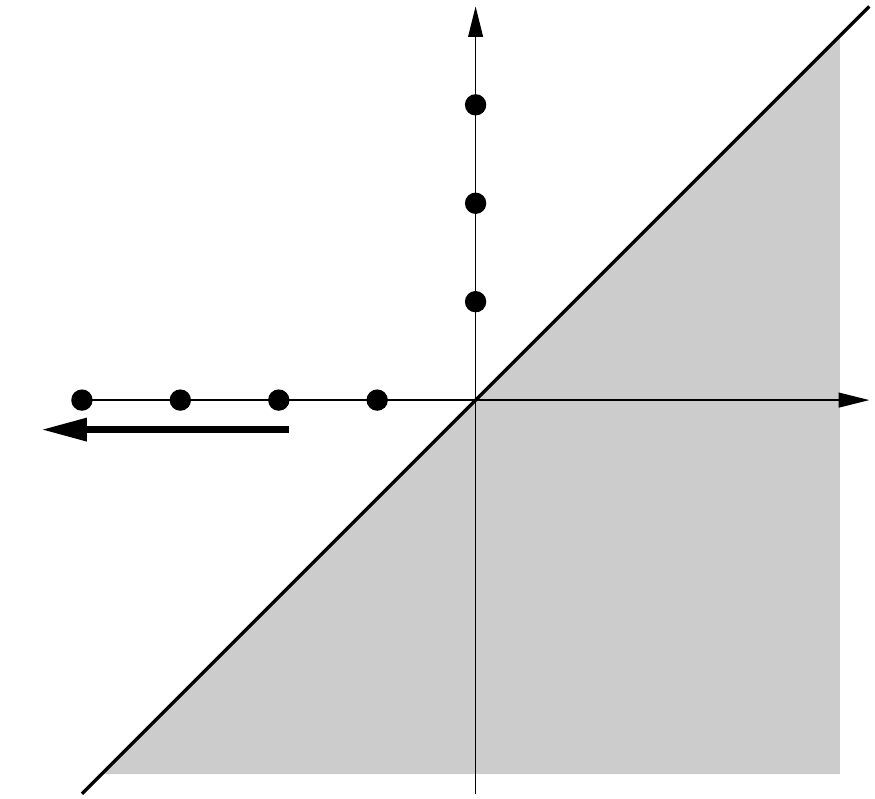}}
  \caption{a) 2d waveguide with bounded inhomogeneity, b) poles of the Laplace transform of the 
    different waveguide modes}
  \end{center}
\end{figure}
For the multi-dimensional case we first present the extension of the
one-dimensional Hardy space method to Helmholtz problems for
waveguides with locally bounded inhomogeneities like the circle in
Fig. \ref{fig:Strip}. Given some incoming wave $u_{\rm i}$ satisfying
the homogeneous Helmholtz equation $-\Delta u_{\rm i} -\kappa^2 u_{\rm
i}=0$, the total wave $u$ is the sum of $u_{\rm i}$ and a scattered
wave $u_{\rm s}$, which has to satisfy a radiating condition. The
problem is given by
\begin{subequations}
\label{eq:Waveguide}
\begin{eqnarray}
  -\Delta u(x,y) - \kappa^2 u(x,y) &=& 0,\quad (x,y) \in \Omega \subset \setR \times [0,\pi] ,\\
  \partial_\nu u &=& 0\qquad \text{on}~\partial \Omega,\\
  \MT \calL u_{\rm s} (\bullet,y) & \in &  H^+(D),\qquad y\in [0,\pi] \label{eq:WGpolecond}.
\end{eqnarray}
\end{subequations}
The Laplace transform in \eqref{eq:WGpolecond} is applied for fixed
tangential variable $y$ to the scattered wave $u_{\rm s}$ outside the
gray shaded domain, i.e. for $x$ in both directions $\pm \infty$.
Note, that $u_{\rm s}$ can be written as a superposition of
one-dimensional waves with wavenumbers $k_n:=\sqrt{\kappa^2 - n^2}$
\begin{equation}
  u_{\rm s}(x,y) = \sum_{n=0}^\infty c_n \cos(n y) e^{i k_n x}.
\end{equation}
As sketched in Fig. \ref{fig:WGpoles} there exist only a few guided
modes and infinitely many evanescent modes. The parameter $\kappa_0$
in the M\"obius transform has to be adapted to the region, where the
Laplace transforms for all these modes are holomorphic. This way
the results of \cite{HohageSchmidtZschiedrich:03a, HohageNannen:08}
carry over to the waveguide problem \eqref{eq:Waveguide}, even though
this is not directly contained in these papers.

In the gray shaded domain $\Omega_{\rm int}:=\Omega \cap [a,b] \times
[0,\pi]$ of Fig. \ref{fig:Strip} we use a standard finite element
method to approximate the total wave $u$. Since only the scattered
wave satisfies the radiation condition, we have to solve the problem
in the exterior domains for $u_{\rm s}$ and not for $u$.  On the
artificial boundary $\Gamma$, that separates $\Omega_{\mathrm int}$
from $\Omega_{\mathrm ext}$ this gives rise to a jump condition for
the Dirichlet values as well as an additional boundary term from the
jump in the Neumann values:
\begin{equation*}
  \begin{aligned}
    &\int_{\Omega_{\rm int}} \left(\nabla u \cdot \nabla v - \kappa^2 u v\right) d(x,y) 
    + \int_0^\pi \int_{-\infty}^a \left(\nabla u_{\rm s} \cdot \nabla v - \kappa^2 u_{\rm s} v\right) dx dy\\
    &+ \int_0^\pi \int_{b}^\infty \left(\nabla u_{\rm s} \cdot \nabla v - \kappa^2 u_{\rm s} v\right) dx dy = 
    \int_0^\pi \left( \left(\partial_x u_{\rm i}\right)(b,y) v(b,y) -  
      \left(\partial_x u_{\rm i}\right)(a,y) v(a,y)\right) dy
  \end{aligned}
\end{equation*}
for suitable test functions $v$. The infinite integrals can be
transformed into the Hardy space using \eqref{eq:basisid} and as in
the one-dimensional case the decomposition \eqref{eq:defOpT}, which
ensures the continuity of the solution over the interfaces. E. g. for
$x\geq b$ the stiffness and mass integral become

\begin{subequations}
  \label{eq:stripintegrals}
  \begin{eqnarray}
    \int_0^\pi \int_{b}^\infty \nabla u_{\rm s} \cdot \nabla v ~ dx dy
    &=&-2 i\kappa_0 \int_b^\infty 
    \BilF{\OpT_+ \vector{{u_{\rm s}}_0(y)}{U(\bullet,y)}}{\OpT_+ \vector{v_0(y)}{V(\bullet,y)}} dy
    \notag 
    \\
    + \frac{2 i }{\kappa_0} \int_b^\infty &&\!\!\!\! \!\!\!\! 
    \BilF{\partial_y \OpT_- \vector{{u_{\rm s}}_0(y)}{U(\bullet,y)}}{\partial_y \OpT_- \vector{v_0(y)}{V(\bullet,y)}} dy,
    \\
    \int_0^\pi \int_{b}^\infty u_{\rm s} ~ v ~ dx dy
    &=& 
    \frac{2 i }{\kappa_0} \int_b^\infty \BilF {\OpT_- \vector{{u_{\rm s}}_0(y)}{U(\bullet,y)}}{\OpT_- \vector{v_0(y)}{V(\bullet,y)}} dy.
  \end{eqnarray}
\end{subequations}
Since on the interface $[0,\pi]$ there already exists a discretization
consisting of the traces $b^{(y)}_m$ of finite element basis functions
in the interior domain, we use these basis functions for the boundary
values ${u_{\rm s}}_0$ and $v_0$:
\begin{equation*}
  {u_{\rm s}}_0 (y) = \sum_{m=0}^{N_m} c_m b^{(y)}_m (y) ,\qquad y \in [0,\pi].
\end{equation*}
For $U,V \in H^+(D) \otimes H^{1/2}([0, \pi])$ it is reasonable to use 
tensor product elements:
\begin{equation*}
  U(z,y) = \sum_{m=0}^{N_m} \sum_{i=0}^{N_i} c_{m,i} z^i b^{(y)}_m (y) ,\qquad z \in S^1,\quad y \in [0,\pi].
\end{equation*}
In this way the integrals~\eqref{eq:stripintegrals} give rise to 
tensor products of the boundary matrices
\begin{equation}
\label{eq.LocBd}
S^{\rm bd}_{mn}=\int_0^\pi \partial_y b^{(y)}_m(y) ~ \partial_y b^{(y)}_n(y) ~ dy, 
\qquad 
M^{\rm bd}_{mn}  = \int_0^\pi  b^{(y)}_m(y) ~  b^{(y)}_n(y) ~ dy
\end{equation}
and the Hardy space matrices $S^{\rm HSM}= -2 i \kappa_0 \OpT_{N,+}^{\top}\OpT_{N,+}$ 
and $M^{\rm HSM} = \frac{2 i}{ \kappa_0} \OpT_{N,-}^{\top}\OpT_{N,-}$:
\begin{equation}
\label{eq.Ext}
S_{\rm ext} = S^{\rm HSM} \otimes M^{\rm bd} + M^{\rm HSM} \otimes S^{\rm bd},\qquad 
M_{\rm ext} = M^{\rm HSM} \otimes M^{\rm bd}.
\end{equation}

\subsection{Tensor product elements}
In the general multi-dimensional case we consider the Helmholtz equation
\begin{equation}
    -\Delta u(x,y) - \kappa^2 n(x,y) u(x,y) = 0,\quad (x,y) \in \Omega,
    \label{eq:Helmholtz} 
\end{equation}
with a potential $n$, an unbounded domain $\Omega$, with some
boundary condition on $\partial \Omega$ and as a radiation condition 
the pole condition along a generalized radial direction for $u$. Soon it 
will become clear what is meant by a generalized radial direction
of $u$ The assumptions on the potential $n$ will be given in Remark
\ref{Rem:poss_potential}.

In \cite{HohageNannen:08} the computational domain is obtained by
intersecting $\Omega$ with a ball $B_a$ of radius $a$, $\Oi:=B_a \cap
\Omega$, such that the unbounded exterior is $\Oe:=\setR^d \setminus
B_a$. Using polar coordinates in $\Oe$ and separation of variables the
unbounded radial direction and the bounded surface directions
separate. Hence, the one-dimensional approach can be applied to the
radial part of the exterior solution and a standard finite element
method handles the interior part as well as the bounded surface
directions. 
\begin{figure}
  \begin{center}
  \subfigure{\resizebox{.4\textwidth}{!}{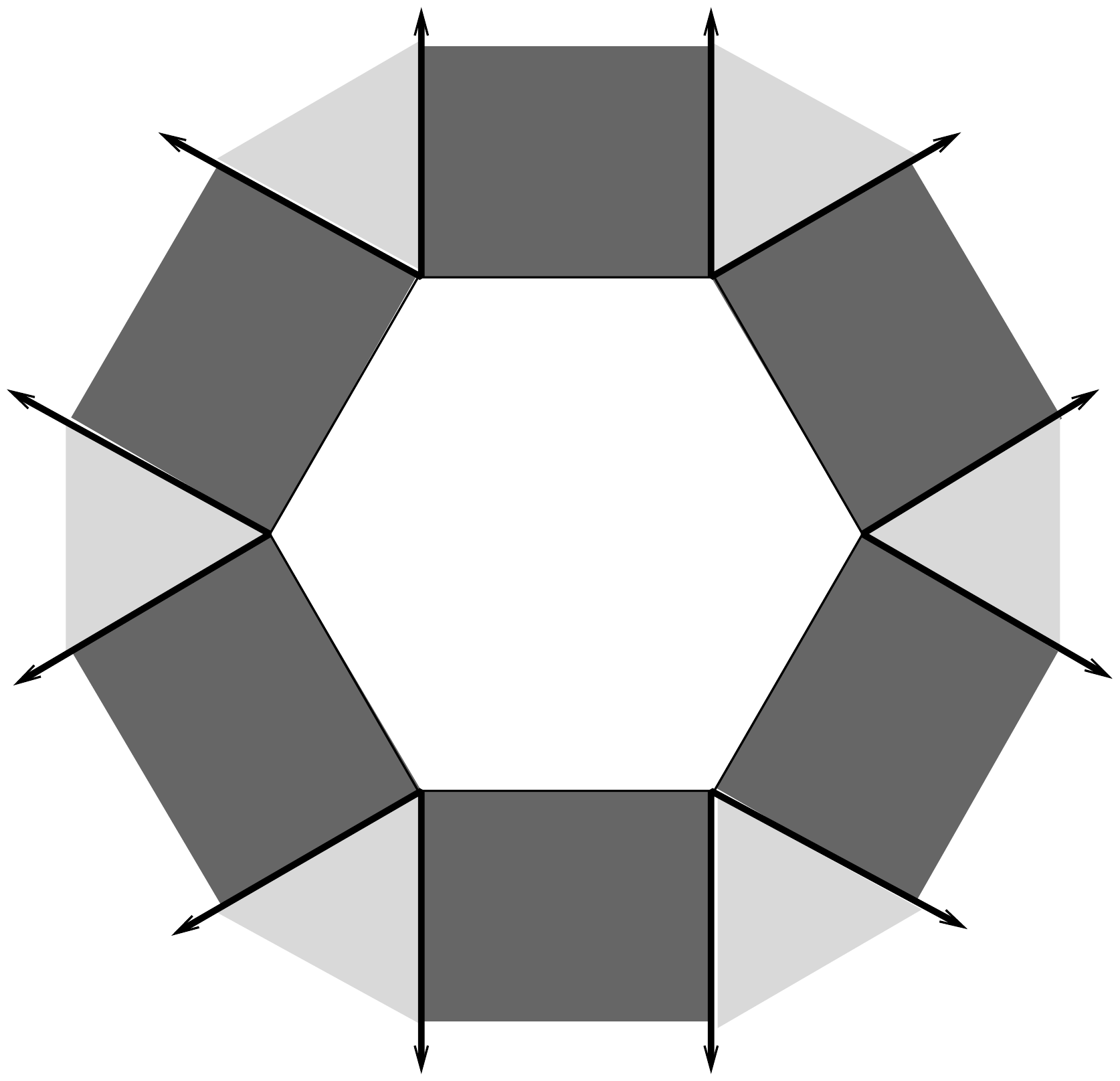}} \hspace{1cm}
  \subfigure{\resizebox{.4\textwidth}{!}{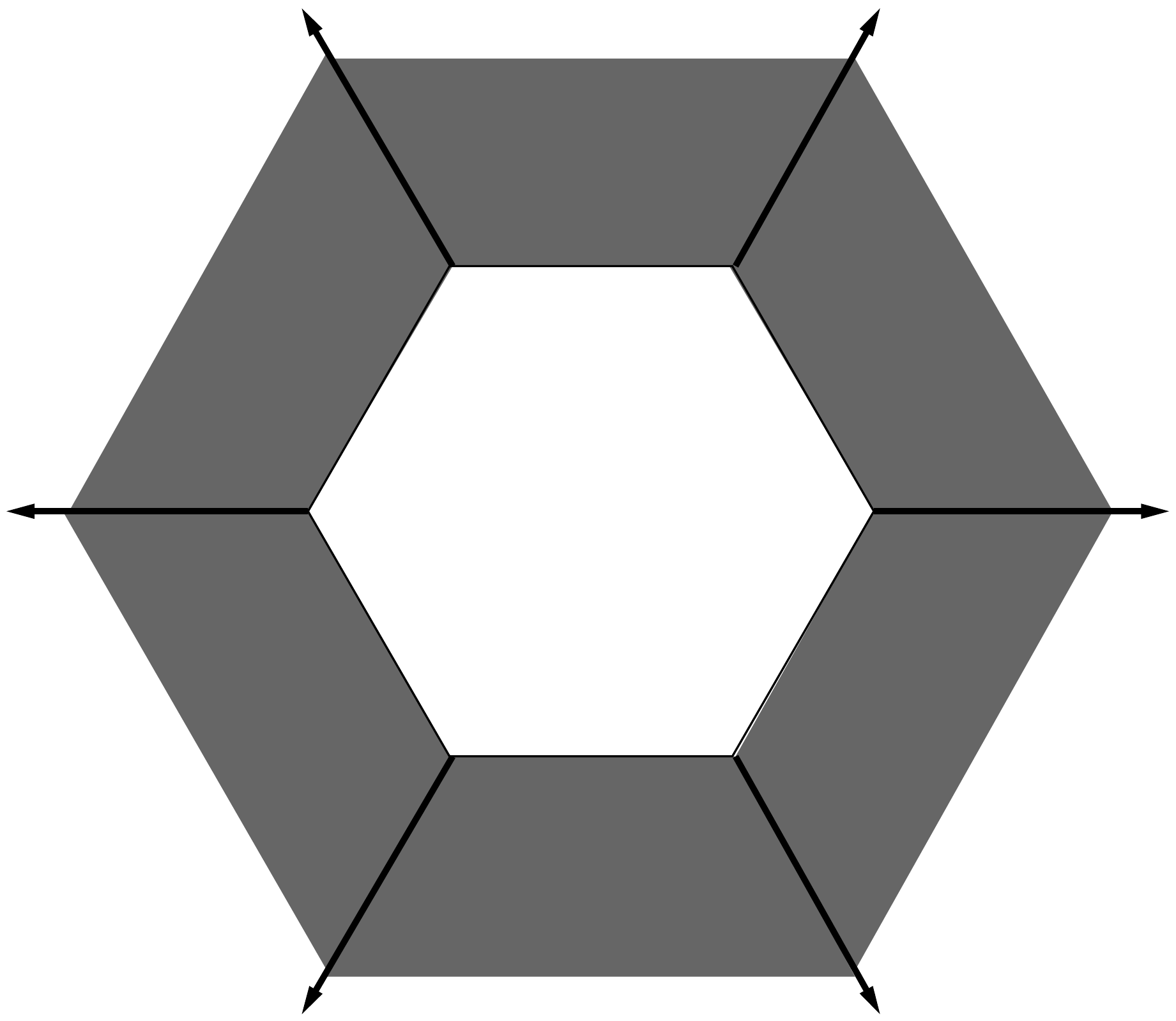}}
  \label{fig:infTriangle}
  \caption{different segmentations of the exterior domain}
  \end{center}
\end{figure}

However the boundary of $\Oe$ need not be a sphere, arbitrary convex
polygons $P$ can be used to split $\Omega$ into $\Oe:=\Omega \setminus
P$ and $\Oi:=P \cap \Omega$ with interface $\Gamma:=\partial P$. Here
for simplicity we present only the two-dimensional case. In
Fig.~\ref{fig:infTriangle} two different segmentations of
the exterior domain are illustrated: The left one decomposes the exterior
into semi-infinite strips and infinite triangles, while the right one uses
semi-infinite trapezoids.

\subsubsection{Infinite triangles and strips}
\label{sec:Triangle}
We first present the Hardy space infinite element method for the
segmentation using infinite triangles and strips. For the
semi-infinite strips we use the elements from
section~\ref{sec:TensorStrip}. Hence only the implementation of the
method for the infinite triangles is considered below. If $P$ is the
vertex of such an infinite triangle and $n_1$ and $n_2$ are the
unit normal vectors of the neighboring strips, then the triangle is
given by $T=\bfg \left([0,\infty) \times [0,\infty)\right)$ with the
linear mapping
\begin{equation}
\bfg (\xi,\eta)= P + \xi n_1 +\eta n_2,\qquad (\xi,\eta) \in [0,\infty) \times [0,\infty).
\end{equation}
If we define $\hat u(\xi,\eta) := u(\bfg(\xi,\eta))$ on the reference, then
with the constant Jacobi matrix $J$ of $\bfg$, the 
mass and stiffness integrals for the infinite triangle are transformed according to
\begin{equation}
  \label{eq:lokalIntTriangle}
  \begin{aligned}
  \int_{T} u ~ v\, d(x,y)  &=   
  \int_{[0,\infty)\times[0,\infty)}  \hat u  ~ \hat v  ~|J| \, d (\xi,\eta),\\
  \int_{T} \nabla_{xy}  u \cdot \nabla_{xy} v \,d(x,y)  &= 
  \int_{[0,\infty)\times[0,\infty)}
  J^{-T} \nabla_{\xi\eta} \hat u \cdot J^{-T} \nabla_{\xi\eta} \hat v |J| \, d (\xi,\eta).
  \end{aligned}
\end{equation}
In contrast to the strips, the integrals include the Jacobi matrix and
there are two infinite directions to which \eqref{eq:basisid} is
applied to. If we define the constant matrix $G:=|J| J^{-1} J^{-T}$,
then the local element matrices for each infinite triangle $T$ are given by
\begin{equation}
  \begin{aligned}
    S_T =& G_{11} S^{\rm HSM} \otimes S^{\rm HSM} +G_{12} S^{\rm HSM} \otimes M^{\rm HSM}\\
    &+ G_{21} M^{\rm HSM} \otimes S^{\rm HSM}+ G_{22} M^{\rm HSM} \otimes M^{\rm HSM}, \\
    M_T =& |J| M^{\rm HSM} \otimes M^{\rm HSM}.
  \end{aligned}
\end{equation}

\begin{Rem}
In the transformation $\bfg$ we use the unit normal vectors, which 
guarantees the continuity of the solution along the infinite rays.
\end{Rem}


\subsubsection{Infinite trapezoids}
\label{sec:Trapezoids}
The method presented in~\ref{sec:Triangle} has the advantage to be
easy to implement, but it results in extra degrees of freedom in the
infinite triangles. The second methods avoids these degrees of freedom
by using infinite trapezoids. The infinite rays, which are no longer
normal to the boundary, could e.g. be constructed in 2d with bisecting
lines. Another possibility is to choose a reference point $P_0$ in the
interior domain and construct the rays $R$ for each vertex $V$ of the
boundary by $R=V-P_0$. General conditions for suitable segmentations
in infinite trapezoids may be found in \cite{Schmidt:02,Zschiedrichetal:06,Kettner:07}.

\begin{figure}
  \resizebox{1\textwidth}{!}{\input{Files/trapez.pdftex_t}}
  \caption{Transformation of each trapezoid}
    \label{fig:transformTrapez}
\end{figure}
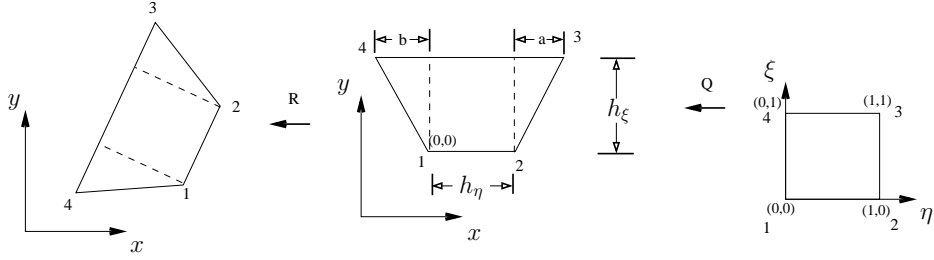
Given a segmentation with finite trapezoids,  
the trapezoid element is the image of a reference rectangle 
(see Fig. \ref{fig:transformTrapez}) under the affine bilinear mapping $\bfg$ with
\begin{equation}
  \label{eq:trafo} (x,y) = \bfg(\eta,\xi)=R \circ Q(\eta,\xi) + (x_1,y_1)^T,
\end{equation}
where 
\begin{equation}
 (x,y)  = Q(\eta,\xi) = \left(
    \begin{array}{c}
      h_\eta \eta- b \xi + (a+b) \eta \xi \\
      h_\xi \xi \\
    \end{array} \right)
\end{equation} 
with $h_\eta = \sqrt{ (x_2-x_1)^2 + (y_2-y_1)^2}$, $a= (x_4-x_3,
y_4-y_3) (x_2-x_3, y_2-y_3)^T /||(x_4-x_3, y_4-y_3)||_2$, $b=
(x_3-x_4, y_3-y_4) (x_1-x_4, y_1-y_4)^T/||(x_4-x_3, y_4-y_3)||_2 $ and
$h_\xi = \sqrt{||(x_3-x_2,y_3-y_2)||_2^2 - a^2}$. Note, that $a$ and
$b$ are signed distance variables.

\begin{Rem}
The variable $\xi$ plays the role of a generalized radial variable,
whereas $\eta$ is the surface variable on $\Gamma$. To guarantee
continuity of the discrete solution in the exterior domain it is
important, that the radial variable $\xi$ along the rays of the
segmentation is independent of the neighboring infinite elements.
This is fulfilled, if two neighboring (finite) trapezoids have the
same boundary vertices.
\end{Rem}

The rotation $R$ is given by
\begin{equation}
  \label{eq:rot}
  (\hat x, \hat y)  = R (\tilde x,\tilde y) =
  \frac{1}{\sqrt{(x_2-x_1)^2+(y_2-y_1)^2}}
  \left(\begin{array}{cc}
      x_2-x_1 & y_1-y_2 \\
      y_2-y_1 & x_2-x_1
    \end{array}
  \right).
\end{equation}
The Jacobi matrix $J$ of the transformation $\bfg$ and its determinant are
\begin{equation}
  \label{eq:Jacobian}
  J= \left(
    \begin{array}{cc}
      h_\eta+(a+b)\xi& -b+(a+b)\eta 
      \\
       0& h_\xi
    \end{array}\right), \qquad
  |J| = h_{\xi}(h_{\eta}+(a+b)\xi) 
\end{equation}
and the inverse of $J$ is 
\begin{equation}
  J^{-1} = \left( 
    \begin{array}{cc}
      \displaystyle \frac{1}{h_{\eta} +
        \xi(a+b)} & \displaystyle \frac{b - (a+b)\eta}{h_{\xi} (h_{\eta} + \xi(a+b))} 
      \\
  0& \displaystyle \frac{1}{h_{\xi}} 
       \end{array}
  \right).
\end{equation}
Note, that $J$ is no longer constant. Mass ans stiffness integral transform as 
in~\eqref{eq:lokalIntTriangle}
\begin{equation}
  \label{eq:lokalInt}
  \begin{aligned}
    \int_{T} u ~ v\, d(x,y)  &=   
    \int_{[0,1]\times[0,\infty]} \hat u ~ \hat v |J| \, d (\eta,\xi),\\
    \int_{T} \nabla_{xy} u \cdot \nabla_{xy} v \,d(x,y)  &= 
     \int_{[0,1]\times[0,\infty]}
    J^{-T} \nabla_{\eta\xi} \hat u \cdot 
    J^{-T} \nabla_{\eta\xi} \hat v ~|J| \, d (\eta,\xi).
  \end{aligned}
\end{equation}
As before we use tensor product basis functions on the reference element 
$\hat b_{m,i}:=\hat b_m^{(\eta)} \otimes \hat b_i^{(\xi)}$. Again the integrals 
decouple, such that for the mass integral we obtain
\begin{equation}
  \label{eq:massintegraltransformed}
  \begin{aligned}
    \int_{T} b_{m,i} ~ b_{n,j} \, d(x,y) =
     \int_0^1 \hat b_m^{(\eta)} \hat b_n^{(\eta)} \, d\eta
      \int_0^\infty \hat b_i^{(\xi)}\hat b_j^{(\xi)}
      h_{\xi}(h_{\eta}+(a+b)\xi) \, d\xi.
  \end{aligned}
\end{equation}
Due to the inverse Jacobian $J^{-1}$ the stiffness integral is more 
complicated. We obtain
\begin{equation}
  \begin{aligned}
    \label{eq:stiffintegraltransformed}
    \int_{T}& \nabla_{xy} b_{m,i} \cdot \nabla_{xy} 
    b_{n,j} \,d(x,y) =  \\
    & \int_0^1 \partial_\eta \hat b_m^{(\eta)} 
      \left(h_\xi^2 + (b-(a+b)\eta)^2 \right)
      \partial_\eta \hat b_n^{(\eta)} \, d\eta 
      \int_0^\infty \frac{\hat b_i^{(\xi)} \hat b_j^{(\xi)}}{h_\xi(h_\eta+(a+b)\xi)} \, 
      d\xi 
    \\
    &+\int_0^1 \partial_\eta \hat b_m^{(\eta)} \frac{b-(a+b)\eta}{h_\xi} \hat b_n^{(\eta)} \, d\eta 
    \int_0^\infty  \hat b_i^{(\xi)} \partial_\xi \hat b_j^{(\xi)} \, d\xi \\
    &+\int_0^1 \hat b_m^{(\eta)} \frac{b-(a+b)\eta}{h_\xi} \partial_\eta \hat b_n^{(\eta)} \, d\eta 
    \int_0^\infty  \partial_\xi \hat b_i^{(\xi)}\hat b_j^{(\xi)} \, d\xi \\ 
    &+ \int_0^1 \hat b_m^{(\eta)} \hat b_n^{(\eta)} \, d\eta 
    \int_0^\infty \partial_\xi \hat b_i^{(\xi)}\frac{h_\eta+(a+b)\xi}{h_\xi} \partial_\xi \hat b_j^{(\xi)} \, d\xi.
\end{aligned}
\end{equation}
For functions depending on $\eta$ we use again the traces of finite
element function in $\Omega_{\rm int}$. Hence, the bounded integrals
over $\eta$ can be treated in the usual way by quadrature
formulas. For the infinite integrals we apply the identity
\eqref{eq:basisid} and the decomposition \eqref{eq:defOpT} to
transform the $\xi$-direction into the Hardy space $H^+(D)$.\par

We have to take care of the factors $\xi$ and $(\xi+c)^{-1}$ with a
constant $c>0$ in \eqref{eq:massintegraltransformed} and
\eqref{eq:stiffintegraltransformed}. For these we need one additional
operator $\Dop:H^+(D) \to H^+(D)$ implicitly defined by
\begin{equation*} 
  \MT \calL \{ (\bullet ) f \}=\MT \left\{ -\left( \calL f\right)' \right\}=\Dop \left(\MT \calL f\right).
\end{equation*} 
Direct calculations yield
\begin{equation}
  \label{eq:defDop} \left(\Dop F \right)(z) =\frac{(z-1)^2}{2 i
    \kappa_0} F'(z)+\frac{z-1}{2 i \kappa_0} F(z),\qquad F \in
  H^+(D).
\end{equation}
If we use the set of trigonometric monomials up to the order $N_\xi$ 
as basis functions in $H^+(D)$, we get the discrete operator
\begin{equation}
  \label{def:diskrDt} \Dop_{N_\xi} := 
  \frac{1}{2 i
    \kappa_0} {\scriptstyle \left(\begin{array} {cccccc} -1 & 1 & & & \\ 1 & -3 & 2
        & & \\ &2 & -5 & 3 & \\ & & \ddots & \ddots & \ddots \\ & & & N_\xi & - 2N_\xi-1
      \end{array}\right)}.
\end{equation}
Obviously it holds for the operator $(\Dop+c~ \id)^{-1}$
\begin{equation}
  \label{eq:defDopinv} \MT \calL \left\{ \frac{1}{\bullet +c} f \right\}=
  (\Dop+c~ \id)^{-1} \left(\MT \calL f\right), \qquad c>0.
\end{equation} 
Note that both operators and the matrices
$\Dop_{N_\xi}$ and $\left(\Dop_{N_\xi}+c ~\id_{(N_\xi+1) \times (N_\xi+1)}\right)^{-1}$ are symmetric.\par

\begin{table*}[tbp] \centering
  \begin{tabular}{Sc|Sc|Sc|Sc} 
    space domain & $\MT \calL$ & def. & implementation 
    \\ 
    \hline 
    $\displaystyle f$ & $\displaystyle \frac{1}{i\kappa_0}
    \OpT_-\left(f_0,F\right)^\top $ & \eqref{eq:defOpT} & \eqref{eq:OpTN} 
    \\ 
    \hline
    $\displaystyle f'$ & $\displaystyle \OpT_+\left(f_0,F\right)^\top$ &
    \eqref{eq:OpTp} & \eqref{eq:OpTN} 
    \\
    \hline 
    $\displaystyle (\bullet) f$ &
    $\displaystyle \Dop \left(\MT \calL f \right)$ & \eqref{eq:defDop} &
    \eqref{def:diskrDt}
    \\
    \hline $\displaystyle \frac{1}{\bullet+c} f,~c>0$ &
    $\displaystyle \left(\Dop+c~\id\right)^{-1} \left(\MT \calL f \right)$ & \eqref{eq:defDopinv} &
    num. inv. of $\left(\Dop_{N}+c~\id\right)$
  \end{tabular}
  \caption{Operators for the HSIE method; for the basic identity see
    \eqref{eq:basisid}}
  \label{tab:operators}
\end{table*}

With the operators given in Table \ref{tab:operators} and the identity
\eqref{eq:basisid} we are able to transform all integrals over $\xi$ into the Hardy
space $H^+(D)$. Using a Galerkin ansatz there with monomial basis
functions, local matrices are obtained. For example the integral 
in~\eqref{eq:massintegraltransformed} yields the local mass matrix
\begin{equation}
  M:=-\frac{2 h_\xi}{i \kappa_0} T_{N_\xi,-}^\top 
  \left(h_\eta~\id_{(N_\xi+1)\times (N_\xi+1)}+(a+b) \Dop_{N_\xi}\right) T_{N_\xi,-}.
\end{equation}
In the same way we can treat the integrals in \eqref{eq:stiffintegraltransformed} and
obtain the matrices
\begin{equation}
  \begin{aligned}
    L_{00}&:=-\frac{2 }{ i \kappa_0 h_\xi} T_{N_\xi,-}^\top \left(h_\eta~\id_{(N_\xi+1)\times (N_\xi+1)}+(a+b) 
      \Dop_{N_\xi}\right)^{-1} T_{N_\xi,-},\\
    L_{01}&:=-2~ T_{N_\xi,-}^\top  T_{N_\xi,+}, \\
    L_{10}&:=-2 ~T_{N_\xi,+}^\top  T_{N_\xi,-},\\
    L_{11}&:=-\frac{2i \kappa_0 }{h_\xi} T_{N_\xi,+}^\top \left(h_\eta~\id_{(N_\xi+1)\times (N_\xi+1)}+(a+b) 
      \Dop_{N_\xi}\right) T_{N_\xi,+}.
  \end{aligned}
\end{equation}

\begin{Rem}
  \label{Rem:poss_potential}
  The method is applicable to a wide range of scalar Helmholtz-type
  problems. Since each segment in $\Oe$ is treated separately, unbounded
  inhomogeneities as e.g. waveguides are possible. Even coefficient functions $n$
  with unbounded support are possible, if there exist a segmentation of $\Oe$ such
  that in each segment the function can be written in sums of terms
  $(\xi+a_j)^{k_j}$ (including negative powers) and functions depending on $\eta$
  \begin{equation}
    \hat n(\xi,\eta) = 
    n \circ \bfg (\xi,\eta)=\sum (\xi +a_j)^{k_j} c_j(\eta), \qquad a_j>0 ~\text{for}~k_j<0.
  \end{equation}
\end{Rem}
\begin{Rem}
The exact statement of the tensor product space is due to the decomposition in
\eqref{eq:defOpT} a little bit complicated. Another problem is, that the infinite
integrals have to be bounded. This can be done by choosing test functions $v$, which
decay fast enough and who are dense in the Hardy space after transformation. The
details can be found in \cite{HohageNannen:08}.
\end{Rem}
\section{Perfectly Matched Layer}
\label{sec:PML}
Exterior complex scaling was introduced by Simon~\cite{Simon:79} to facilitate
the mathematical formulation of boundary conditions for the wave functions in quantum
mechanics.  It is shown by Chew and Weedon~\cite{ChewWeedon:94} that
B\'erenger's~\cite{Berenger:94} PML developed for transient Maxwell's equations may
be interpreted as a complex scaling of the exterior solution.  Thus PML can be
regarded as equivalent to exterior complex scaling.

Starting from a discretization of the exterior as in Fig.~\ref{fig:MCRmesh}, the
generalized radial coordinate is scaled by a constant complex factor $\sigma$, such
that scattered outward radiating waves are damped exponentially.  Using this
discretization no special corner conditions are required.  The exponential damping
justifies to truncate the infinite domain some distance away from the boundary of the
computational domain $\Omega_{int}$ and to impose homogeneous Dirichlet or Neumann
boundary condition.  This way the computational domain is surrounded by a finite
layer.  The convergence of the PML{} method for homogeneous exterior domains is analyzed
by Lassas and Sommersalo~\cite{LassasSomersalo:98} for the scattering problem and by
Kim and Pasciak~\cite{KimPasciak:08} for the resonance problem.

In our numerical experiments we compare the HSIE method with the adaptive PML method
described in~\cite{ZschiedrichBurgeretal:05,SchaedleZschiedrichetal:06}.  
The special feature of this PML
is that the thickness of the layer is chosen adaptively based on an a posteriori
estimate of the error introduced by truncating the layer and taking into account the
discretization error of the interior.  The distribution of the grid points is based
on the observation that inside the PML{} short waves that require a fine grid to be
resolved with a certain accuracy are damped much faster than long waves. Long waves
in turn are well resolved on rather coarse grids.

\section{Convergence test: A strip waveguide}
In this section we compare different discretizations of the exterior
domain for the Hardy space infinite element method applied to
\eqref{eq:Helmholtz} with $\kappa=\frac{2 \pi}{1.5}$, the refraction
index
\begin{equation}
n(x,y)=n(y)=\begin{cases} n_2^2 &y \in (-a,a) \\  n_1^2&y\in \setR \setminus (-a,a) \end{cases}
\end{equation}
and $a=0.0365$, $n_1= 1.45$ and $n_2=3.4$. The incoming wave is given by 
$u_{\rm i} (x,y) = v(y) e^{i \kappa_x x}$,  with
\begin{equation}
 v(y)=
\begin{cases} 
C_1 e^{\sqrt{\kappa_x^2 - n_1^2 \kappa^2}y } &,y \leq -a\\
C_2 e^{-i \sqrt{n_2^2 \kappa^2 - \kappa_x^2} y} + C_3 e^{i \sqrt{n_2^2 \kappa^2 - \kappa_x^2} y}\!\!\!&,y \in (-a,a)\\
C_4 e^{-\sqrt{\kappa_x^2 - n_1^2 \kappa^2} y } &,y \geq a
\end{cases}
\end{equation}
for $\kappa_x>0$ and complex coefficients $C_1$, $C_2$, $C_3$ and $C_4$, which have to ensure 
the continuity of $v$ and $v'$ in $\setR$.
\begin{Rem}
The incoming wave should satisfy the Helmholtz equation. Plugging the
ansatz $u_{\rm i}(x,y)=v(y)e^{i \kappa_x x}$ into \eqref{eq:Helmholtz}
leads to the eigenvalue problem
\begin{equation*}
\left(-\partial_y^2   - \kappa^2 n  \right)  v = - \kappa_x^2 v 
\end{equation*}
for the eigenpair $(\kappa_x^2,v) \in \setR \times H^2(\setR)$. If the
jump in the refraction index $n$ is large enough, such an $\kappa_x
\in (n_1\kappa,n_2\kappa)$ exist and the corresponding eigenfunction
$v$ is exponentially decaying for $y \to \pm \infty $ and oscillating
in $(-a,a)$.
\end{Rem}

\begin{figure}
\centering
\subfigure[\label{fig:WGTestSolField}]{\includegraphics[width=0.30\textwidth]{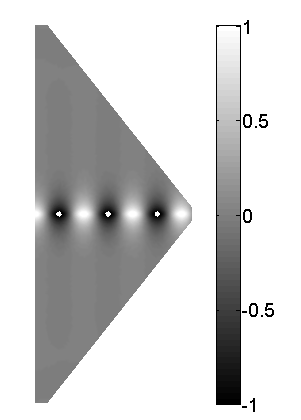}} \hfill
\subfigure[\label{fig:WGTestMesh1}]{\includegraphics[width=0.22\textwidth]{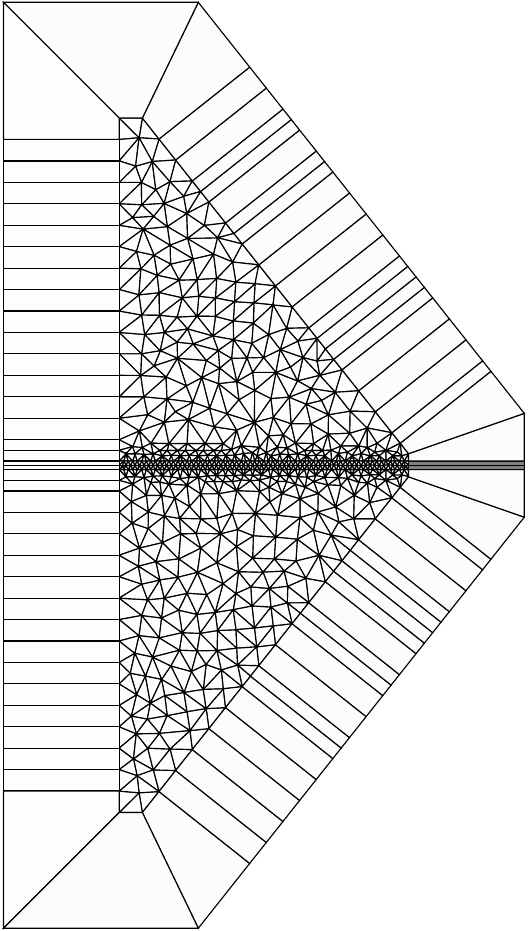}} \hfill
\subfigure[\label{fig:WGTestMesh1a}]{\includegraphics[width=0.22\textwidth]{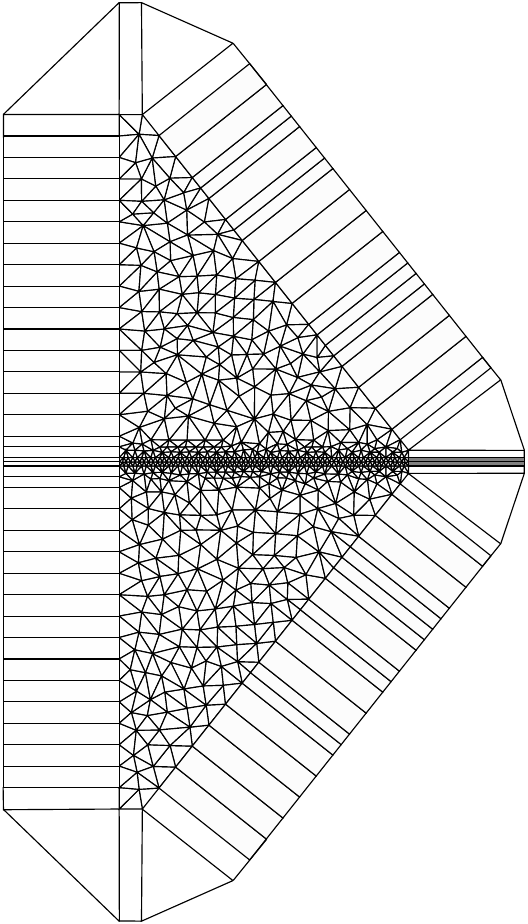}}
\subfigure[\label{fig:WGTestMesh2}]{\includegraphics[width=0.22\textwidth]{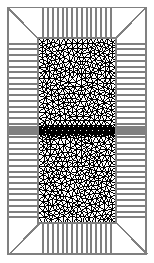}}
\caption{a) real part of $u_{\rm i}$; b,c) 1st mesh with two different exterior discretizations; 
  d) 2nd mesh}
\label{fig:WGTestMesh}
\end{figure}
We solved the problem for strip waveguide for the two different
meshes in Fig. \ref{fig:WGTestMesh}. The incoming wave is coupled via
the jump conditions described in section \ref{sec:TensorStrip} on the
left vertical boundary part of the domain, which is chosen
sufficiently large, so that $u_{\rm i}$ can be set to $0$ in the lower
and upper left corner. In order to test the Hardy space method the
right vertical boundary is for one part of the computations very
small. In the area around the right waveguide port the boundary (and
with it the Hardy space method) has a big influence on the numerical
solution, which should approximate the incoming field $u_{\rm i}$.\par

For the exterior domain two different types of discretization are
used: First we combined infinite strips with infinite triangles
(Sec. \ref{sec:Triangle}). Second we used the trapezoids of section
\ref{sec:Trapezoids}. As shown in Fig. \ref{fig:WGTestMesh} the rays
are chosen in normal direction to the boundary and in corners as
bisecting lines. For the computations we refined the diagrammed meshes
once and used a high order method with polynomial order $7$ for the
finite element method in the interior domain.

\begin{figure}
\centering
\includegraphics[width=\textwidth]{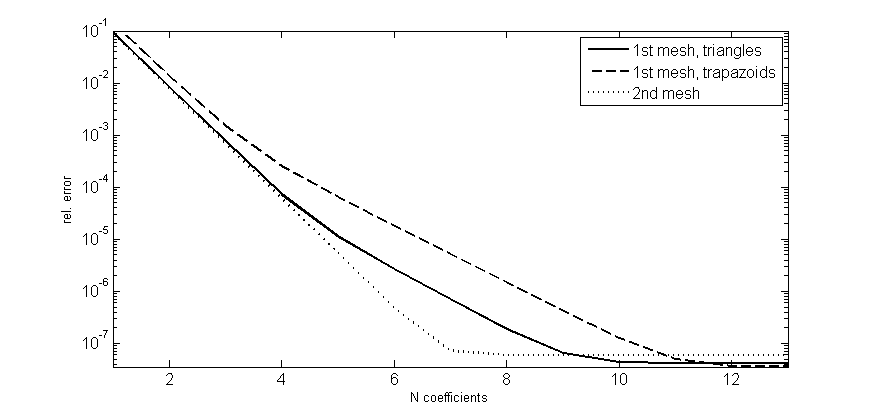}
\caption{Relative $H^1(\Omega_{\rm int})$ error of the Hardy space method vs. 
the number of degrees of freedom in radial direction}
\label{fig:WGTestConverfgence}
\end{figure}

Fig. \ref{fig:WGTestConverfgence} shows for the two different
discretizations and the two meshes in Fig. \ref{fig:WGTestMesh}
exponential convergence of the Hardy space method with respect to the
number of degrees of freedom in radial direction. For the rectangular
mesh both discretizations give the same result, since the influence of
the upper and lower right corners on the discrete solution is very
small. For the challenging mesh the trapezoidal discretization needs
approximately one degree of freedom more than the non-uniform
discretization. The computational costs for both discretizations are
similar.

\begin{figure}
\centering
\subfigure{\includegraphics[width=0.45\textwidth]{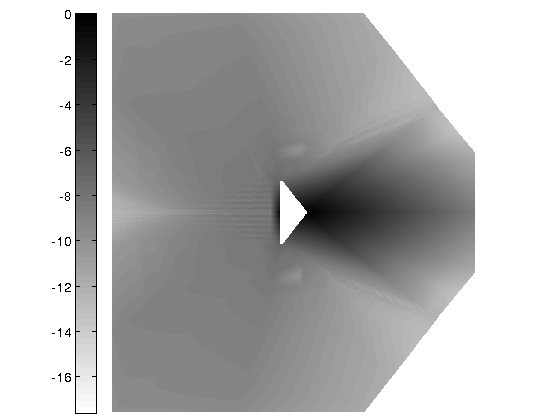}} \hfill
\subfigure{\includegraphics[width=0.45\textwidth]{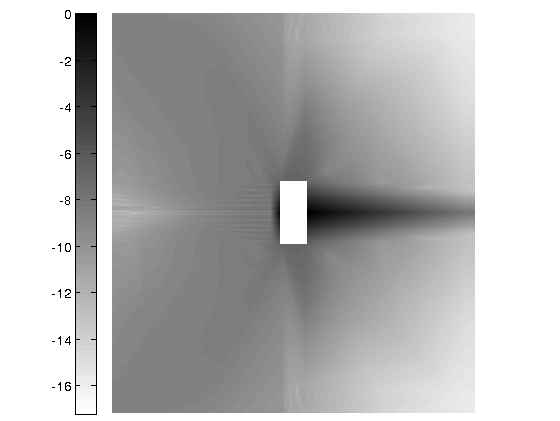}} 
\caption{Logarithm of the absolute value of the degrees of freedom along 
the rays of the exterior domain for two different meshes in Fig. \ref{fig:WGTestMesh}}
\label{fig:WGTestHardyModes}
\end{figure}

Fig. \ref{fig:WGTestHardyModes} shows the logarithm of the absolute
value of the degrees of freedom in the Hardy space. The white domain
in the middle is the interior domain and we have plotted only the
degrees of freedom for the trapezoidal discretization. As expected
from the one-dimensional case and from the theory about spherical
exterior domains they decrease exponentially. Nevertheless, they do
not decrease uniformly in all directions. Especially near to the exit
port of the waveguide the decreasing factor is lower than in the other
regions. For this reason a strategy to choose the number of Hardy
modes adaptively is currently under investigation.

\section{Numerical example: a micro cavity resonator}
\begin{figure}
\centering
\subfigure[\label{fig:MCRsetting}]{\includegraphics[width=0.48\textwidth]{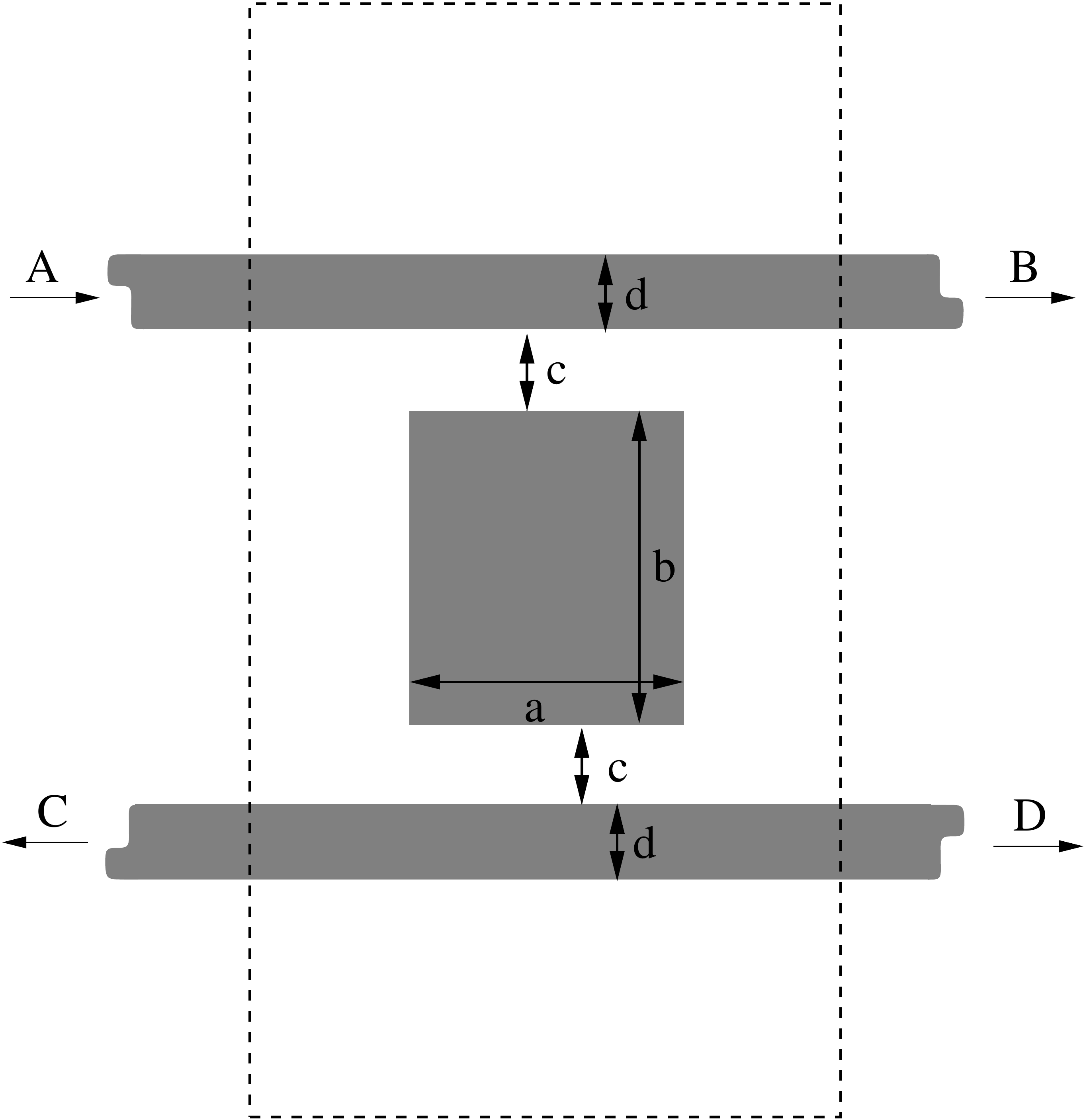}} \hfill
\subfigure[\label{fig:MCRmesh}]{\includegraphics[width=0.44\textwidth]{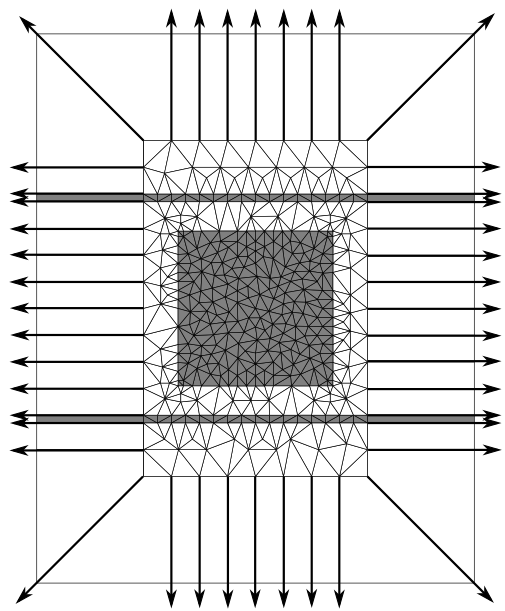}}
\caption{a) Schematic geometry of the micro cavity resonator. b) Coarse grid discretization 
by triangles in the interior and trapezoids in the exterior. 
($a=b=1.451 \mu\mathrm{m}$, 
$c = 0.2745\mu\mathrm{m}$ and $d=0.073\mu\mathrm{m}$)}
\label{fig:MCR}
\end{figure}

This example is taken from Hammer \cite{Hammer:02} and consists of two
waveguides coupled by a cavity. In order to exclude
discretization effects originating from the resolution of the layout a geometry is
chosen that is defined by polygons. The computational domain has a size of $3.5\mu\mathrm{m}
\times 4.546 \mu\mathrm{m}$, with a square cavity $a=b=1.451 \mu\mathrm{m}$ in the center. 
We have $c = 0.2745\mu\mathrm{m}$ and $d=0.073 \mu\mathrm{m}$ in Fig.~\ref{fig:MCRsetting}.


\subsection{Scattering problem}
\label{subsec:scattering}

\begin{figure}
  \centering
  \subfigure{\includegraphics[width=0.40\textwidth]{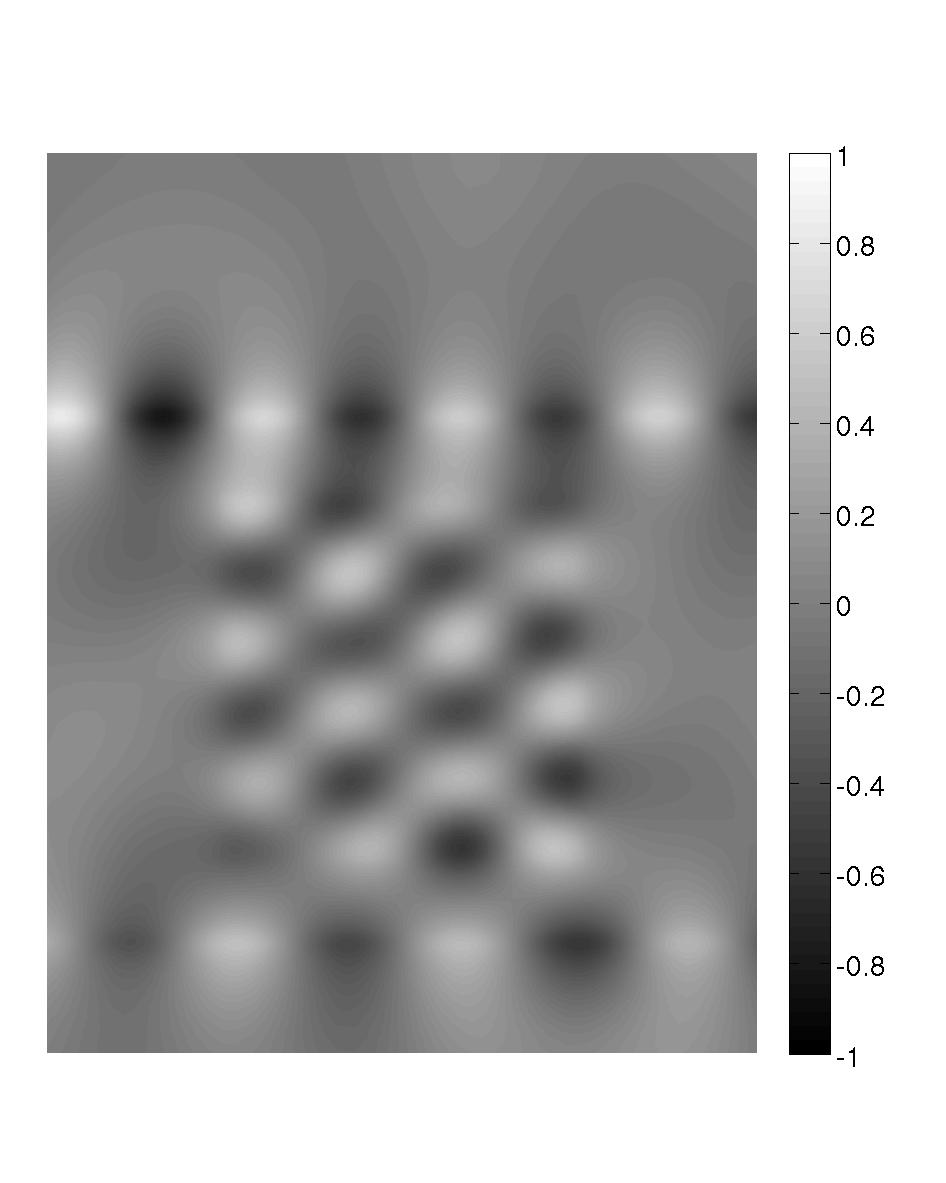}}\hfill
  \subfigure{\includegraphics[width=0.40\textwidth]{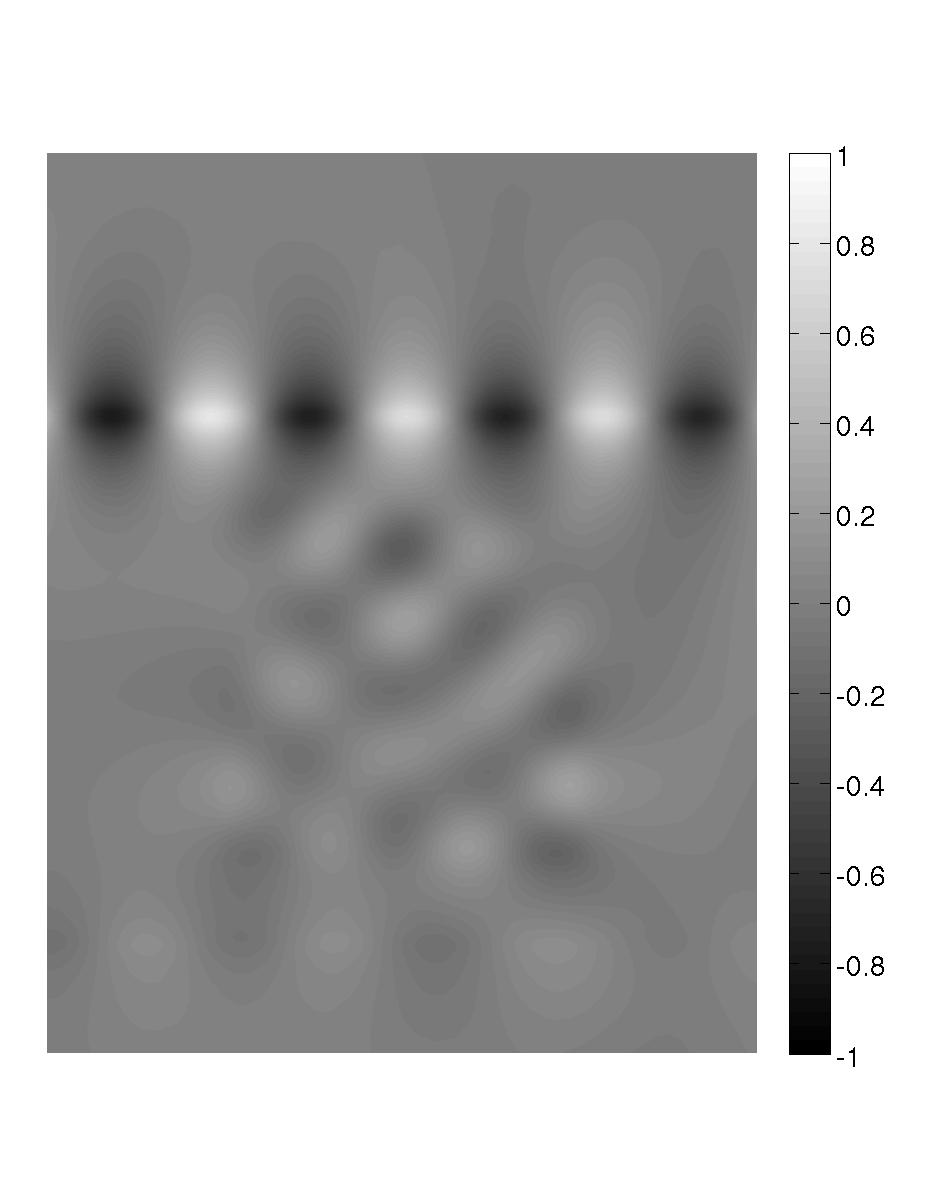}}
  \caption{Real part of a the solution to the scattering problem with wavelength 
    $k=\frac{2 \pi}{1.5759 \mu\mathrm{m}}$ ({\it left}) and 
    $k=\frac{2 \pi}{1.5 \mu\mathrm{m}}$ ({\it right}) \label{fig:solMCR}} 
\end{figure}

To model the scattering of an incoming wave by an object the incoming
waveguide mode of the last section is coupled by a jump condition of
the Neumann and Dirichlet data at the left vertical boundary part for
two different wavenumbers $k_{1}=\frac{2 \pi}{1.5\mu\mathrm{m}}$ and
$k_{2}=\frac{2 \pi}{1.5759 \mu\mathrm{m}}$. $k_{2}$ is close to a
resonance, whereas for $k_{1}$ the cavity has only little effect,
which can be seen in Fig. \ref{fig:solMCR}. For $k_2$ the wave
propagates through the cavity into the lower waveguide.

For these computations we used for the interior finite elements of
$5$th order and refined the coarse grid~\ref{fig:MCRmesh} three times
uniformly. In the exterior domain we used the trapezoidal Hardy space
infinite element method with the parameter $\kappa_0=8+5 i$ and $30$
degrees of freedom in order to discretize the Hardy space. In total we
got approximately 200.000 degrees of freedom.\par

%

%

\subsection{Resonance problem} 
\label{subsec:resonance}

Starting from a coarse mesh the interior is adaptively refined using a 
residual based error estimator~\cite{ZschiedrichBurgeretal:05}. On each refinement 
level the eigenvalue of the resonance problem is calculated.  
To evaluate the relative error in the eigenvalue, a reference resonance
frequency if $\omega=1.1951173e+15 - 1.489202e+13i$ is calculated
on a very fine mesh using the PML{}. The wavelength of the resonance is 
$\lambda = 2\pi c /\omega$, where $c=299792458 m/s$ is the speed of light.
Hence the wavenumber $k_{1} = 2\pi/\lambda_{1} = \frac{2 \pi}{1.5759 \mu\mathrm{m}} $
of Section~\ref{subsec:scattering} is close to the resonance.

In this example two methods, PML and HSIE, to realize transparent 
boundary conditions are compared. 
The results for the PML{} are obtained using
the adaptive PML as described in Section~\ref{sec:PML}. The results for the HSIE 
are obtained using $N=2,\dots,50$ Hardy modes and selecting the best result. 
\begin{figure}
  \centering
  \includegraphics[width=0.49\textwidth]{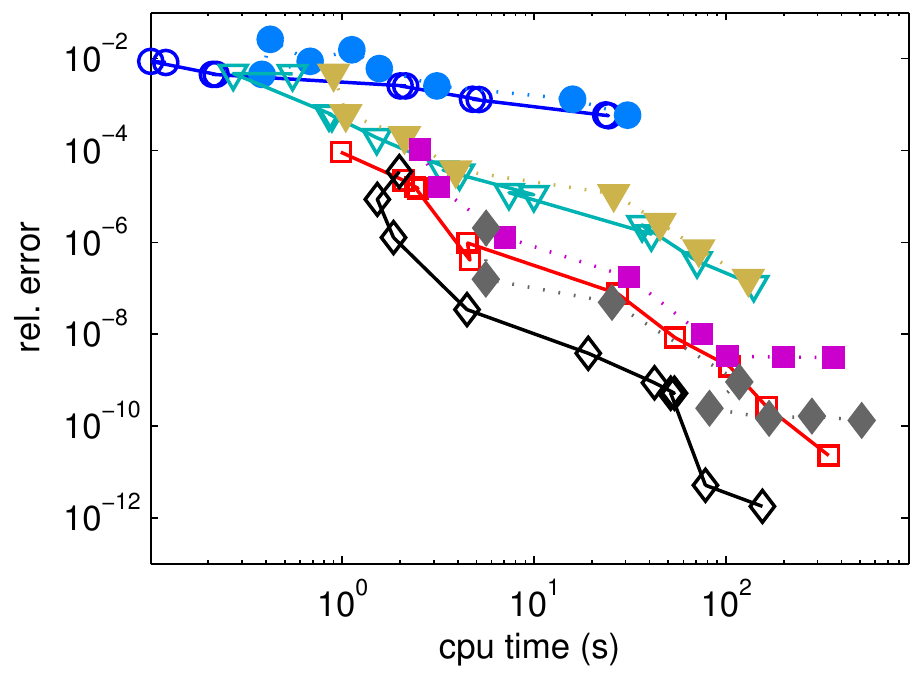}
  \includegraphics[width=0.49\textwidth]{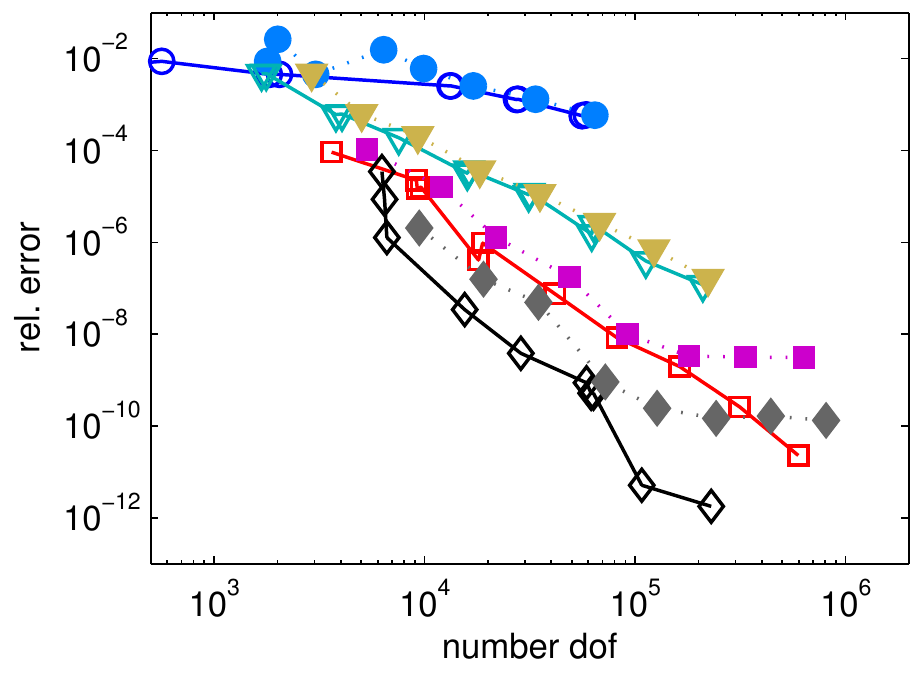}
  \caption{Comparison: PML (dashed line with filled markers) 
    with HSIE (solid line) with $\kappa_{0} = 5+3i$ 
    for various finite element degrees:
    $1$ ($\circ$), $2$ ($\triangle$), $3$ ($\Box$), and $4$ ($\Diamond$).
    Left: Work-precision diagram showing the relative error in the eigenvalue
    vs.  cpu-time in seconds.
    Right: Convergence of PML and HSIE showing the relative error in the 
    eigenvalue vs. number of degrees of freedom.}
  \label{fig:errorcpundof}
\end{figure}
\ref{fig:errorcpundof} shows the error versus the cpu time (left) and the total
number of degrees of freedom (right) that is required to solve the eigenvalue problem
on the finest refinement level.  The initial guess is $1.195\cdot
10^{15}-0.01489\cdot 10^{15}i$. The HSIE method in this example yields results that are better or
at least as good the results obtained with the PML method. 
The implementation was done in the C++ code
JCMsuite~\cite{JCMwave}.

\begin{figure}
  \centering
  \includegraphics[width=0.49\textwidth]{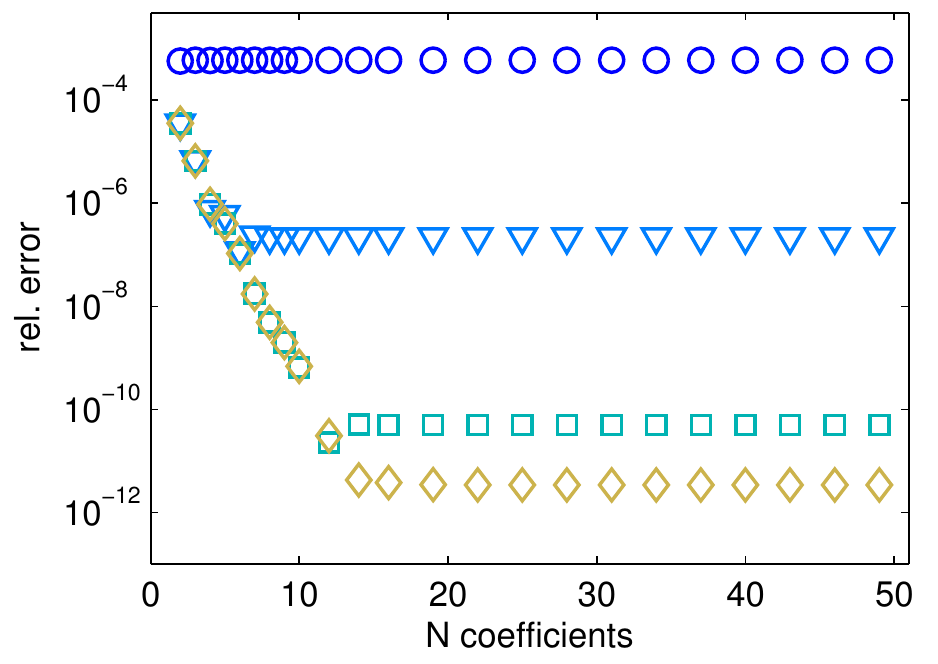}
  \caption{Rel. error in the eigenvalue vs. number of Hardy-modes 
    ($\kappa_{0} = 5+3i$) for finite elements of degree 
    $p=1$ ($\circ$), $p=2$ ($\triangle$), $p=3$ ($\Box$), and $p=4$ ($\Diamond$). 
  }
  \label{fig:errornmodes}
\end{figure}
Fig~\ref{fig:errornmodes} shows the super algebraic convergence in the number of Hardy
modes. For these calculations $\kappa_{0} = 5+3i$.  The number of degrees of freedom to
discretize the interior is about $944000$ for $p=1$, $169000$ for $p=2$, $121000$ for
$p=3$ and $103000$ for $p=4$.  For low accuracies very few Hardy modes
are required, e.g. to reduce the approximation error of the transparent boundary
condition below $10^{-6}$ $7$ Hardy modes are sufficient. To obtain a relative error
of about $10^{-10}$ almost one third of the total number of degrees of freedom is
spent on approximating the transparent boundary condition.

\section{Conclusions}
We have presented here the Hardy space infinite element method from a practical point
of view. In this form inhomogeneous exterior domains can be treated as well as
arbitrary convex polygons as artificial boundaries. The numerical results show
superalgebraic convergence with respect to the degrees of freedom in the Hardy space,
i.e. the degrees of freedom in the radial direction. The method preserves the
eigenvalue structure and is therefore well suited for solving resonance
problems. Compared to former realizations of the pole condition like the {\it cut
function approach} (\cite{Schmidtetal:07}) it is not possible to recover the solution
in the exterior domain directly. In this point the HSIE method behaves like perfectly
matched layer methods: The degrees of freedom in the exterior domain are physically
irrelevant. They just provide for a good approximation at the exact transparent
boundary condition on the artificial boundary.\par A comparison of the results of the
HSIE and the used PML is already given in the last section. From the point of
implementation the HSIE requires a new ({\it in})finite element, whereas the
PML only changes the variational formulation of the problem. On the other hand the
HSIEs show exponentially convergence, while the PML inherits the convergence
order of the used finite element method. One basic difference between both methods is
given by the nature of the discretizations. For the PML there exist two steps:
First truncating the infinite PML-domain and than discretizing the finite layer
using finite elements. The HSIE uses a (transformed) variational formulation of
the whole infinite domain. The only discretization results from the Galerkin
method.\par The HSIE is not restricted to scalar, time-harmonic problems. In
\cite{Ruprechtetal:08} a similar version of the method is used for solving
time-depending problems. Moreover, there exist first results for Maxwell's equations.

\section*{Acknowledgment}
The fruitful and stimulating discussions with T. Hohage, F. Schmidt and
L. Zschiedrich are sincerely acknowledged. 
\bibliographystyle{elsarticle-num} \bibliography{bibliography}
\end{document}

%% file: Files/poles.pdf_t
\begin{picture}(0,0)%
\includegraphics{Files/poles.pdf}%
\end{picture}%
\setlength{\unitlength}{4144sp}%
\begingroup\makeatletter\ifx\SetFigFontNFSS\undefined%
\gdef\SetFigFontNFSS#1#2#3#4#5{%
  \reset@font\fontsize{#1}{#2pt}%
  \fontfamily{#3}\fontseries{#4}\fontshape{#5}%
  \selectfont}%
\fi\endgroup%
\begin{picture}(10114,4734)(339,-3358)
\put(2521,-2401){\makebox(0,0)[lb]{\smash{{\SetFigFontNFSS{20}{24.0}{\rmdefault}{\mddefault}{\updefault}{\color[rgb]{0,0,0}$- i \kappa$}%
}}}}
\put(3151,-601){\makebox(0,0)[lb]{\smash{{\SetFigFontNFSS{20}{24.0}{\rmdefault}{\mddefault}{\updefault}{\color[rgb]{0,0,0}$i \kappa$}%
}}}}
\put(3241,-2851){\makebox(0,0)[lb]{\smash{{\SetFigFontNFSS{20}{24.0}{\rmdefault}{\mddefault}{\updefault}{\color[rgb]{0,0,0}$P_{\kappa_0}^-$}%
}}}}
\put(4276,-1951){\makebox(0,0)[lb]{\smash{{\SetFigFontNFSS{20}{24.0}{\rmdefault}{\mddefault}{\updefault}{\color[rgb]{0,0,0}$\kappa_0$}%
}}}}
\put(5716,1109){\makebox(0,0)[lb]{\smash{{\SetFigFontNFSS{20}{24.0}{\rmdefault}{\mddefault}{\updefault}{\color[rgb]{0,0,0}$m_{\kappa_{0}}(z):=i \kappa_0~ \frac{z+1}{z-1}$}%
}}}}
\put(10261,-1681){\makebox(0,0)[lb]{\smash{{\SetFigFontNFSS{20}{24.0}{\rmdefault}{\mddefault}{\updefault}{\color[rgb]{0,0,0}$1$}%
}}}}
\end{picture}%

%% file: Files/polesWaveguide2.pdf_t
\begin{picture}(0,0)%
\includegraphics{Files/polesWaveguide2.pdf}%
\end{picture}%
\setlength{\unitlength}{4144sp}%
\begingroup\makeatletter\ifx\SetFigFontNFSS\undefined%
\gdef\SetFigFontNFSS#1#2#3#4#5{%
  \reset@font\fontsize{#1}{#2pt}%
  \fontfamily{#3}\fontseries{#4}\fontshape{#5}%
  \selectfont}%
\fi\endgroup%
\begin{picture}(3997,3644)(526,-3233)
\put(800,-1816){\makebox(0,0)[lb]{\smash{{\SetFigFontNFSS{20}{24.0}{\rmdefault}{\mddefault}{\updefault}{\color[rgb]{0,0,0}$n \to \infty$}%
}}}}
\put(2020,-1681){\makebox(0,0)[lb]{\smash{{\SetFigFontNFSS{20}{24.0}{\rmdefault}{\mddefault}{\updefault}{\color[rgb]{0,0,0}$\mathfrak{i} k_3$}%
}}}}
\put(2300,-1096){\makebox(0,0)[lb]{\smash{{\SetFigFontNFSS{20}{24.0}{\rmdefault}{\mddefault}{\updefault}{\color[rgb]{0,0,0}$\mathfrak{i} k_2$}%
}}}}
\put(2300,-646){\makebox(0,0)[lb]{\smash{{\SetFigFontNFSS{20}{24.0}{\rmdefault}{\mddefault}{\updefault}{\color[rgb]{0,0,0}$\mathfrak{i} k_1$}%
}}}}
\put(2300,-196){\makebox(0,0)[lb]{\smash{{\SetFigFontNFSS{20}{24.0}{\rmdefault}{\mddefault}{\updefault}{\color[rgb]{0,0,0}$\mathfrak{i} k_0$}%
}}}}
\end{picture}%

%% file: Files/InfiniteTriangles.pdf_t
\begin{picture}(0,0)%
\includegraphics{Files/InfiniteTriangles.pdf}%
\end{picture}%
\setlength{\unitlength}{4144sp}%
\begingroup\makeatletter\ifx\SetFigFontNFSS\undefined%
\gdef\SetFigFontNFSS#1#2#3#4#5{%
  \reset@font\fontsize{#1}{#2pt}%
  \fontfamily{#3}\fontseries{#4}\fontshape{#5}%
  \selectfont}%
\fi\endgroup%
\begin{picture}(7648,7378)(857,-6495)
\put(3871,-2896){\makebox(0,0)[lb]{\smash{{\SetFigFontNFSS{34}{40.8}{\rmdefault}{\mddefault}{\updefault}{\color[rgb]{0,0,0}$\Omega_{\rm int}$}%
}}}}
\end{picture}%

%% file: Files/infiniteTrapazoids.pdf_t
\begin{picture}(0,0)%
\includegraphics{Files/infiniteTrapazoids.pdf}%
\end{picture}%
\setlength{\unitlength}{4144sp}%
\begingroup\makeatletter\ifx\SetFigFontNFSS\undefined%
\gdef\SetFigFontNFSS#1#2#3#4#5{%
  \reset@font\fontsize{#1}{#2pt}%
  \fontfamily{#3}\fontseries{#4}\fontshape{#5}%
  \selectfont}%
\fi\endgroup%
\begin{picture}(8413,7288)(497,-6405)
\put(3871,-2896){\makebox(0,0)[lb]{\smash{{\SetFigFontNFSS{34}{40.8}{\rmdefault}{\mddefault}{\updefault}{\color[rgb]{0,0,0}$\Omega_{\rm int}$}%
}}}}
\end{picture}%

%% file: Files/trapez.pdftex_t
\begin{picture}(0,0)%
\includegraphics{Files/trapez.pdftex}%
\end{picture}%
\setlength{\unitlength}{4144sp}%
\begingroup\makeatletter\ifx\SetFigFont\undefined%
\gdef\SetFigFont#1#2#3#4#5{%
  \reset@font\fontsize{#1}{#2pt}%
  \fontfamily{#3}\fontseries{#4}\fontshape{#5}%
  \selectfont}%
\fi\endgroup%
\begin{picture}(8304,2410)(361,-1884)
\put(8665,-1555){\makebox(0,0)[lb]{\smash{{\SetFigFont{14}{12.0}{\rmdefault}{\mddefault}{\updefault}{\color[rgb]{0,0,0}$\eta$}%
}}}}
\put(7247,-259){\makebox(0,0)[lb]{\smash{{\SetFigFont{14}{12.0}{\rmdefault}{\mddefault}{\updefault}{\color[rgb]{0,0,0}$\xi$}%
}}}}
\put(5816,-614){\makebox(0,0)[lb]{\smash{{\SetFigFont{14}{12.0}{\rmdefault}{\mddefault}{\updefault}{\color[rgb]{0,0,0}$h_\xi$}%
}}}}
\put(4470,-1311){\makebox(0,0)[lb]{\smash{{\SetFigFont{14}{12.0}{\rmdefault}{\mddefault}{\updefault}{\color[rgb]{0,0,0}$h_\eta$}%
}}}}
\put(4544,-1739){\makebox(0,0)[lb]{\smash{{\SetFigFont{14}{12.0}{\rmdefault}{\mddefault}{\updefault}{\color[rgb]{0,0,0}$x$}%
}}}}
\put(3357,-381){\makebox(0,0)[lb]{\smash{{\SetFigFont{14}{12.0}{\rmdefault}{\mddefault}{\updefault}{\color[rgb]{0,0,0}$y$}%
}}}}
\put(1474,-1849){\makebox(0,0)[lb]{\smash{{\SetFigFont{14}{12.0}{\rmdefault}{\mddefault}{\updefault}{\color[rgb]{0,0,0}$x$}%
}}}}
\put(361,-516){\makebox(0,0)[lb]{\smash{{\SetFigFont{14}{12.0}{\rmdefault}{\mddefault}{\updefault}{\color[rgb]{0,0,0}$y$}%
}}}}
\end{picture}%

%% file: HIEvsPML.bbl
\begin{thebibliography}{10}

\bibitem{Astley:00}
R.~J. Astley.
\newblock Infinite elements for wave problems: A review of current formulations
  and an assessment of accuracy.
\newblock {\em Internat. J. Numer. Methods Engrg.}, 49(7):951--976, 2000.

\bibitem{Berenger:94}
J.-P. Berenger.
\newblock {A perfectly matched layer for the absorption of electromagnetic
  waves.}
\newblock {\em J. Comput. Phys.}, 114(2):185--200, 1994.

\bibitem{ChewWeedon:94}
W.~C. Chew and W.~H. Weedon.
\newblock A 3d perfectly matched medium from modified {M}axwell's equations
  with stretched coordinates.
\newblock {\em Microwave Optical Tech. Letters}, 7:590--604, 1994.

\bibitem{DemkowiczGerdes:98}
L.~Demkowicz and K.~Gerdes.
\newblock Convergence of the infinite element methods for the {H}elmholtz
  equation in separable domains.
\newblock {\em Numer. Math.}, 79:11--42, 1998.

\bibitem{Duren:70}
P.~L. Duren.
\newblock {\em Theory of {$H\sp{p}$} spaces}.
\newblock Pure and Applied Mathematics, Vol. 38. Academic Press, New York,
  1970.

\bibitem{Givoli:04}
D.~Givoli.
\newblock High-order local non-reflecting boundary conditions: a review.
\newblock {\em Wave Motion}, 39:319--326, 2004.

\bibitem{GroteKeller:98}
M.~J. Grote and J.~B. Keller.
\newblock Nonreflecting boundary conditions for {M}axwell's equation.
\newblock {\em J. Comput. Phys.}, 139:327--324, 1998.

\bibitem{Hammer:02}
M.~Hammer.
\newblock Resonant coupling of dielectric optical waveguides via rectangular
  microcavities: the coupled guided mode perspective.
\newblock {\em Optics Communications}, 214(1-6):155--170, 2002.

\bibitem{HeinHohageKochSchoeberl:07}
S.~Hein, T.~Hohage, W.~Koch, and J.~Sch\"oberl.
\newblock Acoustic resonances in high lift configuration.
\newblock {\em J. Fluid Mech.}, 582:179--202, 2007.

\bibitem{Hoffman:62}
K.~Hoffman.
\newblock {\em Banach spaces of analytic functions}.
\newblock Prentice-Hall Series in Modern Analysis. Prentice-Hall Inc.,
  Englewood Cliffs, N. J., 1962.

\bibitem{HohageNannen:08}
T.~Hohage and L.~Nannen.
\newblock Hardy space infinite elements for scattering and resonance problems.
\newblock {\em SIAM J. Num. Analysis}, 47(2):972--996, 2009.

\bibitem{Hohageetal:02}
T.~Hohage, F.~Schmidt, and L.~Zschiedrich.
\newblock A new method for the solution of scattering problems.
\newblock In B.~Michielsen and F.~Decav\`ele, editors, {\em Proceedings of the
  JEE'02 Symposium}, pages 251--256, Toulouse, 2002. ONERA.

\bibitem{HohageSchmidtZschiedrich:03a}
T.~Hohage, F.~Schmidt, and L.~Zschiedrich.
\newblock Solving time-harmonic scattering problems based on the pole
  condition. {I}. {T}heory.
\newblock {\em SIAM J. Math. Anal.}, 35(1):183--210 (electronic), 2003.

\bibitem{HsiaoWendland:08}
G.~C. Hsiao and W.~L. Wendland.
\newblock {\em Boundary integral equations}, volume 164 of {\em Applied
  Mathematical Sciences}.
\newblock Springer-Verlag, Berlin, 2008.

\bibitem{JCMwave}
{JCMwave GmbH}.
\newblock www.jcmwave.com.

\bibitem{Kettner:07}
B.~Kettner.
\newblock {E}in {A}lgorithmus zur prismatoidalen {D}iskretisierung von
  unbeschr\"ankten {A}u\ss en\-r\"aumen in 2{D} und 3{D}.
\newblock Master's thesis, Freie Universit\"at Berlin, 2007.

\bibitem{KimPasciak:08}
S.~Kim and J.~E. Pasciak.
\newblock The computation of resonances in open systems using a perfectly
  matched layer.
\newblock {\em Math. Comp.}, 2008.

\bibitem{LassasSomersalo:98}
M.~Lassas and E.~Somersalo.
\newblock On the existence and the convergence of the solution of the pml
  equations.
\newblock {\em Computing}, 60:229--241, 1998.

\bibitem{LenoirVulliermeHazard:92}
M.~Lenoir, M.~Vullierme-Ledard, and C.~Hazard.
\newblock Variational formulations for the determination of resonant states in
  scattering problems.
\newblock {\em SIAM J. Math. Anal.}, 23:579--608, 1992.

\bibitem{Nannen:08}
L.~Nannen.
\newblock {\em {H}ardy{-R}aum {M}ethoden zur numerischen {L\"o}sung von
  {S}treu- und {R}esonanzproblemen auf unbeschr{\"a}nkten {G}ebieten}.
\newblock PhD thesis, University of G{\"o}ttingen, {D}er {A}ndere {V}erlag,
  T{\"o}nning, 2008.

\bibitem{Ruprechtetal:08}
D.~Ruprecht, A.~Sch{\"a}dle, F.~Schmidt, and L.~Zschiedrich.
\newblock Transparent boundary conditions for time-dependent problems.
\newblock {\em SIAM J. Sci. Comput.}, 30(5):2358--2385, 2008.

\bibitem{SchaedleZschiedrichetal:06}
A.~Sch\"adle, L.~Zschiedrich, S.~Burger, R.~Klose, and F.~Schmidt.
\newblock Domain decomposition method for {M}axwell's equations: Scattering off
  periodic structures.
\newblock {\em J. Comput. Phys.}, 226:477--493, 2007.

\bibitem{Schmidt:98}
F.~Schmidt.
\newblock An alternative derivation of the exact dtn-map on a circle.
\newblock Technical Report SC 98-32, Konrad-Zuse-Zentrum Berlin, 1998.

\bibitem{Schmidt:02}
F.~Schmidt.
\newblock A new approach to coupled interior-exterior {H}elmholtz-type
  problems: Theory and algorithms.
\newblock Habilitation, Freie Universität Berlin, 2002.

\bibitem{Schmidtetal:07}
F.~Schmidt, T.~Hohage, R.~Klose, A.~Sch{\"a}dle, and L.~Zschiedrich.
\newblock Pole condition: A numerical method for {H}elmholtz-type scattering
  problems with inhomogeneous exterior domain.
\newblock {\em J. Comput. Appl. Math.}, 218(1):61--69, 2008.

\bibitem{Simon:79}
B.~Simon.
\newblock The definition of molecular resonance curves by the method of
  exterior complex scaling.
\newblock {\em Phys. Lett. A}, 71A(2, 3), 1979.

\bibitem{ZschiedrichBurgeretal:05}
L.~Zschiedrich, S.~Burger, R.~Klose, A.~Sch\"adle, and F.~Schmidt.
\newblock {JCMmode}: An adaptive finite element solver for the computation of
  leaky modes.
\newblock In Y.~Sidorin and C.~A. W\"achter, editors, {\em Integrated Optics
  IX}, volume 5728 of {\em Proc. SPIE}, pages 192--202, 2005.

\bibitem{Zschiedrichetal:06}
L.~Zschiedrich, R.~Klose, A.~Sch{\"a}dle, and F.~Schmidt.
\newblock A new finite element realization of the perfectly matched layer
  method for {H}elmholtz scattering problems on polygonal domains in two
  dimensions.
\newblock {\em J. Comput. Appl. Math.}, 188(1):12--32, 2006.

\end{thebibliography}
